\documentclass[12pt,fleqn,leqno]{article}
\usepackage{amsmath,amssymb}

\usepackage[dvips]{graphicx}
\usepackage[dvips]{color}
\usepackage{fancyhdr}

\allowdisplaybreaks

\numberwithin{equation}{section}

\makeatletter
\renewcommand{\section}{\@startsection{section}{1}{0pt}{20pt}{6pt}{\large\bf}}
\renewcommand{\@seccntformat}[1]{\csname the#1\endcsname.\ }

\def\footnoterule{\kern -3pt \hrule width 2.7 true cm \kern 2.6pt}

\def\h{\hspace}
\def\np{\newpage}
\def\nt{\noindent}
\def\p{\!+\!}
\def\m{\!-\!}

\def\EE{\mathsf E\:\!}
\def\PP{\mathsf P}

\def\R{I\!\!R}
\def\LL{I\!\!L}
\def\eps{\varepsilon}

\begin{document}

\title{\bf Global $C^1$ Regularity of the Value Function \\ in Optimal
Stopping Problems}
\author{T.\ De Angelis\, \&\, G.\ Peskir}
\date{}
\maketitle






{\par \leftskip=2cm \rightskip=2cm \footnotesize

We show that if either the process is strong Feller and the boundary
point is probabilistically regular for the stopping set, or the
process is strong Markov and the boundary point is probabilistically
regular for the interior of the stopping set, then the boundary
point is Green regular for the stopping set. Combining this
implication with the existence of a continuously differentiable flow
of the process we show that the value function is continuously
differentiable at the optimal stopping boundary whenever the gain
function is so. The derived fact holds both in the parabolic and
elliptic case of the boundary value problem under the sole
hypothesis of probabilistic regularity of the optimal stopping
boundary, thus improving upon known analytic results in the PDE
literature, and establishing the fact for the first time in the case
of integro-differential equations. The method of proof is purely
probabilistic and conceptually simple. Examples of application
include the first known probabilistic proof of the fact that the
time derivative of the value function in the American put problem is
continuous across the optimal stopping boundary.

\par}





\vspace{-20pt}

\section{Introduction}

A challenging question in boundary value problems is to establish
\emph{regularity} of the solution up to the boundary. By regularity
we mean continuity, differentiability, and/or higher degrees of
smoothness. The problem has a long and venerable history. Continuity
results can be traced back to Poincar\'e \cite{Poi} and the
references therein. Differentiability results date back to Gevrey
\cite{Ge} for parabolic equations and Kellog \cite{Ke-1} for
elliptic equations (see also \cite{Ke-2}). Extensions to more
general parabolic and elliptic equations were made possible using
the techniques developed by Schauder \cite{Sc} (see \cite{Ko} for
further details). As a rule of thumb in the PDE literature it is
known that (probabilistic) regularity of the boundary implies
continuity of the solution up to the boundary, and smoothness (or
H\"older continuity) of the boundary implies smoothness of the
solution up to the boundary (see e.g.\ \cite[Theorem 7, p.\
64]{Fr-1} for parabolic equations and \cite[Lemma 6.18, p.\ 111]{GT}
for elliptic equations). This common belief translates to free
boundary problems for parabolic and elliptic equations as well (see
e.g.\ \cite[Lemma 4.5, p.\ 167]{Fr-2} for a definite result of this
kind dating back to Gevrey \cite{Ge} as well as \cite{Ca} and
\cite[Chapter 8]{CS} for rela- ted results in higher dimensions).
The analytic method of variational inequalities removes the focus
from the free boundary itself and derives a global continuity of the
\emph{space} derivative (for parabolic and elliptic equations of
\emph{diffusion} processes) when the obstacle function is
\emph{globally} $C^1$ while establishing that the time derivative
exists in a \emph{weak} sense only (see \cite[Corollary 1.3, p.\
207]{BL} and \cite[Theorem 3.2, p.\ 26; Theorem 8.2, p.\ 77; Theorem
8.4, p.\ 80]{Fr-3}). The latter fact is not surprising since the
time derivative can fail to exist in the absence of probabilistic
regularity of the free boundary (see e.g.\ \cite[Example 14]{Pe-4}).
A probabilistic approach in \cite{Ok} returns to a probabilistic
regularity of the free boundary by assuming moreover that the free
boundary is twice continuously differentiable and thus making the
assumption `intractable' as the paper points out itself.

In this paper we develop a conceptually simple/direct probabilistic
method which shows that the differentiability results for free
boundary problems can be derived solely from a probabilistic
regularity of the boundary i.e.\ with no need for its smoothness (or
H\"older continuity) of any kind. This applies to (i) both the space
derivative and the time derivative, (ii) more general strong
Markov/Feller processes (not just diffusions), and (iii) both smooth
and non-smooth obstacle functions. Free boundary problems (in
analysis) are known to be equivalent to optimal stopping problems
(in probability) and we derive the differentiability results in the
context of optimal stopping problems which are also of interest in
themselves. We do that by establishing a \emph{continuous} smooth
fit between the value function and the gain (obstacle) function at
the optimal stopping (free) boundary that is traditionally derived
using probabilistic methods in a directional sense only (see Section
2 for details).

In Section 2 we formulate the optimal stopping problem
\eqref{2.1}/\eqref{2.2} and explain its background in terms of (i)
strong Markov/Feller processes, (ii) boundary point regularity
(probabilistic, Green, barrier, Dirichlet), (iii) stochastic flow
regularity, and (iv) infinitesimal generator regularity (including
continuous and smooth fit). In Section 3 we show that if either the
process is strong Feller and the boundary point is probabilistically
regular for the stopping set, or the process is strong Markov and
the boundary point is probabilistically regular for the interior of
the stopping set, then the boundary point is Green regular for the
stopping set (in the sense that the expected waiting time for
entering the stopping set vanishes as the initial point of the
process approaches the boundary point from within the continuation
set). Combining this implication with the existence of a
continuously differentiable flow of the process we show in Sections
4 and 5 that the value function is continuously differentiable at
the optimal stopping boundary whenever the gain function is so.
Theorems 8 and 10 deal with the space derivative (in infinite and
finite horizon respectively) and Theorems 13 and 15 deal with the
time derivative (in infinite and finite horizon respectively).
Examples 12 and 17 derive the analogous regularity results for the
space derivative and the time derivative respectively, when the gain
function is not smooth away from the optimal stopping boundary,
using the local time of the process on the singular points at which
the smoothness breaks down.

The advantage of the probabilistic method employed in the derived
results is that the only hypothesis on the optimal stopping boundary
used is its probabilistic regularity for the stopping set or its
interior (which is implied by monotonicity of the optimal stopping
boundary for instance). This level of generality is insufficient for
the PDE methods as they require at least a Lipschitz (or H\"older)
continuity of the optimal stopping boundary. The derived results
hold both in the parabolic and elliptic case of the free boundary
problem, thus improving upon known analytic results in the PDE
literature, and establishing the fact for the first time in the case
of integro-differential equations. Moreover, the `lifting' method of
Example 17 to our knowledge is applied for the first time in the
literature. It enables one to `lift' a Lipschitz continuity of the
superharmonic/value function to its $C^1$ regularity at Green
regular boundary points. Among other implications this yields the
first known probabilistic proof of the fact that the time derivative
of the value function in the American put problem is continuous
across the optimal stopping boundary.

In parallel to producing a first draft of the present paper we have
also applied/tested some parts of the method of proof in specific
examples. This includes \cite{DeAn} for the time derivative in the
Brownian motion case and \cite{JP} for the space derivative in the
Bessel process case. For further/existing applications to (singular)
stochastic control problems and optimal stopping games we refer to
\cite{DE} and \cite{DGV} respectively. Among intermediate references
we note that the paper \cite{BDM} studies continuity of the time
derivative of solutions to parabolic free-boundary problems in one
(spatial) dimension under the hypotheses that $G=0$ on the stopping
set (with $G>0$ at the end of time) and $H<0$ globally in the
optimal stopping problem \eqref{2.2} below. These hypotheses are
rarely satisfied in the mainstream examples of optimal stopping
problems studied in the literature (including the American put
problem where $G>0$ on the stopping set and $H=0$ globally) and the
present paper fills this gap as well.

\section{Problem formulation}

In this section we introduce the setting of the problem and explain
its background in terms of the general hypotheses imposed and
sufficient conditions that imply them.

\vspace{6pt}

1.\ \emph{Optimal stopping problem}. We consider the optimal
stopping problem
\begin{equation} \h{3pc} \label{2.1}
V(x) = \sup_\tau\, \EE_x \Big[\;\! e^{-\Lambda_\tau} G(X_\tau)
+ \int_0^\tau e^{-\Lambda_t} H(X_t)\, dt\;\! \Big]
\end{equation}
for $x \in \R^d$ with $d \ge 1$ where $X=(X^1, \ldots ,X^d)$ is a
standard Markov process (in the sense of \cite[p.\ 45]{BG}) taking
values in $\R^d$. Thus $X$ is strong Markov, right-continuous with
left limits, and left-continuous over stopping times. The process
$X$ starts at $x$ under the probability measure $\PP_{\!x}$ for $x
\in \R^d$ (or its measurable subset identified with $\R^d$ in the
sequel for simplicity). The supremum in \eqref{2.1} is taken over
all stopping times $\tau$ of $X$ (i.e.\ stopping times with respect
to the natural filtration of $X$), or equivalently, over all
stopping times $\tau$ with respect to a (right-continuous)
filtration $({\cal F}_t)_{t \ge 0}$ that makes $X$ a strong Markov
process under $\PP_{\!x}$ for $x \in \R^d$. All stopping times
considered throughout are assumed to be finite valued unless
otherwise stated (upon recalling that extensions to infinite valued
stopping times are both standard and straightforward). We will also
consider the optimal stopping problem \eqref{2.1} with finite
horizon obtained by imposing an upper bound $T>0$ on $\tau$. In this
case we also need to account for the length of the remaining time so
that \eqref{2.1} extends as follows
\begin{equation} \h{3pc} \label{2.2}
V(t,x) = \sup_{0 \le \tau \le T-t} \, \EE_x \Big[\;\! e^{-\Lambda_\tau}
G(X_\tau) + \int_0^\tau e^{-\Lambda_s} H(X_s)\, ds\;\! \Big]
\end{equation}
for $t \in [0,T]$ and $x \in \R^d$. Note that this includes the case
when the functions $G$ and $H$ are time dependent which can be
formally obtained by setting $X_t^1 = t$ for $t \ge 0$. The
functional $\Lambda$ in \eqref{2.1} and \eqref{2.2} is defined by
\begin{equation} \h{8.5pc} \label{2.3}
\Lambda_t = \int_0^t \lambda(X_s)\, ds
\end{equation}
where $\lambda$ is a continuous function with values in
$[0,\infty)$. The real-valued functions $G$ and $H$ are also assumed
to be continuous. Under these hypotheses it is known (cf.\ \cite{PS}
and \cite{Sh}) that the first entry time of $X$ into the (finely)
closed set $D$ where $V$ equals $G$ (the stopping set) is optimal in
\eqref{2.1}/\eqref{2.2} provided that $G(X)$ and $H(X)$ satisfy mild
integrability conditions. This is true for example if $\lambda>0$
and both $G$ and $H$ are bounded but this sufficient condition can
be considerably strengthened (see \cite{PS} and \cite{Sh} for
details). The (finely) open set where $V$ is strictly larger than
$G$ (the continuation set) will be denoted by $C$. The (optimal
stopping) boundary between the sets $C$ and $D$ will be denoted by
$\partial C$. We will make use of and distinguish between the
\emph{first entry time} of $X$ into $D$ defined by
\begin{equation} \h{8pc} \label{2.4}
\tau_D = \inf\, \{\, t \ge 0\; \vert\; X_t \in D\, \}
\end{equation}
and the \emph{first hitting time} of $X$ to $D$ defined by
\begin{equation} \h{8pc} \label{2.5}
\sigma_D = \inf\, \{\, t > 0\; \vert\; X_t \in D\, \}
\end{equation}
where $D$ can also be replaced by any other measurable subset of
$\R^d$ and an upper bound applies to admissible $t$ in \eqref{2.4}
and \eqref{2.5} when the horizon is finite as in \eqref{2.2}. When
the standard regularity hypotheses recalled above are satisfied, or
any other sufficient conditions implying that $\tau_D$ is optimal in
\eqref{2.1}/\eqref{2.2}, we will say that the problem
\eqref{2.1}/\eqref{2.2} is \emph{well posed}. This will be a
standing premise for the rest of the paper. Any additional
hypotheses will always be invoked explicitly in the statements of
the results below when needed.

\vspace{6pt}

2.\ \emph{Strong Feller processes}. Recall that the process $X$ is
\emph{strong Feller} if
\begin{equation} \h{8pc} \label{2.6}
x \mapsto \EE_x \big[ F(X_t) \big]\;\; \text{is continuous}
\end{equation}
for every real-valued (bounded) \emph{measurable} function $F$ with
$t>0$ given and fixed. Recall also that $X$ is \emph{Feller} if
\eqref{2.6} holds for every real-valued (bounded) \emph{continuous}
function $F$. Recall finally that Feller processes are strong
Markov. Strong Feller processes were introduced and initially
studied by Girsanov \cite{Gi}. All one-dimensional diffusions $X$ in
the sense of It\^o and McKean \cite{IM} are known to be strong
Feller processes because the transition density $p$ of $X$ with
respect to its speed measure $m$ (in the sense that $\PP_{x}(X_t\!
\in\! dy) = p(t;x,y)\, m(dy)$\!) can be chosen to be jointly
continuous in all three arguments (cf.\ \cite[p.\ 149]{IM}). Unique
weak solutions to (non-degenerate) SDEs driven by a Wiener process
in $\R^d$ are known to be not only strong Markov but also strong
Feller processes (see e.g.\ \cite[p.\ 170]{RW}). A time-space
process such as $((t,W_t))_{t \ge 0}$ where $W$ is a standard Wiener
process is \emph{not} a strong Feller process. Not all L\'evy
processes are strong Feller either. Hawkes \cite[Theorem 2.2]{Ha}
showed that a L\'evy process $X$ is strong Feller if and only if
$\PP_{\!x}(X_t\! \in\! dy) \ll \ell(dy)$ for every $t>0$ and $x \in
\R^d$ where $\ell$ denotes Lebesgue measure on $\R^d$. Strong Feller
property is important in relation to boundary point regularity. We
will now present basic facts in this direction.

\vspace{6pt}

3.\ \emph{Boundary point regularity}. There are four closely related
concepts of boundary point regularity that we will address in the
sequel. Throughout we let $b(c,r)$ denote the open ball in the
Euclidean topology of $\R^d$ with centre at $c$ and radius $r>0$. By
$\bar C$ we denote the closure of $C$ and by $D^\circ$ we denote the
interior of $D$. Recall that a real-valued function $v$ is
superharmonic on a set $A \subseteq \R^d$ relative to $X$ if
$\EE_x[v(X_\tau)] \le v(x)$ for all $x \in A$ and all (bounded)
stopping times $\tau \le \tau_{A^c}$ of $X$. A boundary point $z \in
\partial C$ is said to be:
\begin{align} \h{0pc} \label{2.7}
&\text{\emph{Probabilistically regular} (PR) if we have}\;\; \PP_{
\!z}(\sigma_D = 0) = 1\, ; \\[6pt] \label{2.8} &\text{\emph{Green
regular} (GR) if we have}\! \lim_{\;\;C \ni x \rightarrow z \in
\partial C} \PP_{\!x}(\tau_D \ge \eps) = 0\;\; \text{for each}\;\;
\eps>0\, ; \\[2pt] \label{2.9} & \text{\emph{Barrier regular} (BR)
if there exists a superharmonic function}\;\; v>0\;\; \text{on}
\\[2pt] \notag &b(z,r) \cap C\;\; \text{relative to}\;\; X\;\;
\text{for some}\;\; r>0\;\; \text{such that}\! \lim_{\;\;C \ni
x \rightarrow z \in \partial C} v(x) = 0\, ; \\[2pt] \label{2.10}
&\text{\emph{Dirichlet regular} (DR) if}\! \lim_{\;\; C \ni x
\rightarrow z \in \partial C} \EE_x[F(X_{\tau_D})] = F(z)\;\;
\text{for each real-valued (boun-}\\[-4pt] \notag &\text{ded)
measurable function}\;\; F\;\; \text{on}\;\; b(z,r) \cap \bar
C\;\; \text{with}\;\; r>0 \;\; \text{that is continuous at}
\;\; z\, .
\end{align}

\nt Regularity of $z \in \partial C$ in definitions
\eqref{2.7}-\eqref{2.10} refers to the set $D$. If we replace $D$ in
\eqref{2.7}-\eqref{2.10} by any measurable subset $A$ of $\R^d$ then
we speak about regularity of $z \in \partial C$ for the set $A$. By
Blumenthal's 0-1 law (cf.\ \cite[p.\ 30]{BG}) we know that the
probability in \eqref{2.7} can only be either zero or one. The
super(harmonic) function $v$ in \eqref{2.9} is referred to as a
\emph{barrier} itself. The main example of a barrier is $v(x) =
\EE_x(\tau_D)$ for $x \in b(z,r) \cap C$ with $r>0$ when $\lim_{\,C
\ni x \rightarrow z \in \partial C} v(x) = 0$ holds (where $\tau_D$
could be replaced by $\tau_D \wedge 1$ to make it bounded).

It is well known (cf.\ \cite[pp.\ 32-40]{Dy-2}) that if $X$ is
strong Feller then
\begin{equation} \h{8pc} \label{2.11}
\text{PR} \Longleftrightarrow \text{GR} \Longleftrightarrow \text{BR}\, .
\end{equation}
Moreover, if $X$ is strong Feller and uniformly continuous on
compacts in the sense that
\begin{equation} \h{7pc} \label{2.12}
\lim_{t \downarrow 0}\, \sup_{x \in K}\, \PP_{\!x} \Big( \sup_{s \in [0,t]}
\vert X_s \m x \vert > \eps \Big) = 0
\end{equation}
for each compact set $K$ in $\R^d$ and each $\eps>0$ then
\begin{equation} \h{7pc} \label{2.13}
\text{PR} \Longleftrightarrow \text{GR} \Longleftrightarrow \text{BR}
\Longleftrightarrow \text{DR}
\end{equation}
where $\vert \cdot \vert$ denotes the Euclidean norm in $\R^d$. We
will see in the proofs below that our main focus will be on the
\emph{Green regularity}. When the process $X$ fails to be strong
Feller however, then the first equivalence in
\eqref{2.11}/\eqref{2.13} can break down generally, and we will then
require probabilistic regularity for $D^\circ$ instead of $D$ to
gain the Green regularity. Further details in this direction will be
presented in the next section.

We will close this subsection with a few historical details aimed at
clarifying definitions \eqref{2.7}-\eqref{2.10} above. Note that
many papers cited below contain \emph{sufficient conditions} for
boundary point regularity that are directly relevant for the main
results in Sections 4 and 5 below.

Definition \eqref{2.7} embodies what probabilists understand under
regularity. Definition \eqref{2.10} embodies what analysts
understand under regularity. The implication
\eqref{2.7}$\Rightarrow$\eqref{2.10} was first proved by Doob
\cite{Do} for a Wiener process and was then extended by Girsanov
\cite{Gi} to other strong Feller processes. The converse implication
\eqref{2.10}$\Rightarrow$\eqref{2.7} for strong Feller processes was
derived by Krylov \cite {Kr}. Definition \eqref{2.8} embodies a
``hybrid'' condition representing a mixture of \eqref{2.7} and
\eqref{2.10} that makes it suitable for applications as we will see
below. Definition \eqref{2.9} is often used to derive various
sufficient conditions for regularity. Poincar\'e \cite{Poi} used
barriers to derive a sphere condition. Zaremba \cite{Za} replaced
sphere by a cone (cf.\ \cite[pp.\ 247-250]{KS}). Wiener \cite{Wi}
derived a necessary and sufficient condition for regularity using
the capacity of a set (Wiener's test). These papers deal with the
Laplace equation (when $X$ is a Wiener process) and extensions to
more general elliptic equations are normally not difficult
(probabilistically this can be seen through time changes and
comparison arguments). The same phenomenon does not hold for the
heat equation (when $X$ is a time-space Wiener process) and more
general parabolic equations (see e.g.\ \cite[Theorem 8.1]{EK} for a
simple example). Petrovsky \cite{Pe} derived sufficient conditions
for regularity in the heat equation by considering boundaries as
functions of time (Kolmogorov-Petrovsky's test). Necessary and
sufficient conditions for regularity in the heat equation were
announced by Landis \cite{Land}. An analogue of Wiener's test for
the heat equation was derived in the papers by Lanconelli
\cite{Lanc} and Evans \& Gariepy \cite{EG} (see pp.\ 295-296 in the
latter paper for related results and historical comments). We refer
to the paper by Watson \cite{Wa} and the references therein for
subsequent analytic results and further developments. Boundary point
regularity and continuity of the solution to the Dirichlet problem
for standard Markov processes have been studied by Dembinski
\cite{Dem} using purely probabilistic methods (see also the
references therein for further probabilistic papers on this topic).

\vspace{6pt}

4.\ \emph{Stochastic flow regularity}. Stochastic processes whose
sample paths are indexed by their initial points are referred to as
\emph{stochastic flows}. Motivated by needs in the proofs below we
will assume that the standard Markov process $X$ can be realised as
a stochastic flow $(X_t^x)_{t \ge 0, x \in \R^d}$ on a probability
space $(\Omega,{\cal F},\PP)$ in the sense that $\text{Law}(X\,
\vert\, \PP_{\!x}) = \text{Law}(X^x\, \vert\, \PP)$ where we set
$X^x = (X_t^x)_{t \ge 0}$ for $x \in \R^d$.

Examples of stochastic flows include a standard Wiener flow $W =
(W_t^x)_{t \ge 0, x \in \R}$ where $W_t^x = x \p W_t$ (which extends
to all L\'evy processes analogously), an exponential Wiener flow $S
= (S_t^x)_{t \ge 0, x \in \R}$ where $S_t^x = x \exp (\sigma W_t \p
(\mu \m \sigma^2\!/2) t)$ for $\sigma>0$ and $\mu \in \R$, and a
reflecting Wiener flow $R = (R_t^x)_{t \ge 0, x \in \R}$ where
$R_t^x = x \vee \sup_{\,0 \le s \le t} W_s \m W_t$. Very often an
explicit construction of the stochastic flow is not possible and
then one usually aims to establish its existence satisfying
some/further regularity properties. Among these we will need to
consider \emph{continuous}, \emph{differentiable}, and
\emph{continuously differentiable} stochastic flows. For us in this
paper it will mean that there exists a (universal) set $N \in \cal
F$ satisfying $\PP(N)=0$ such that the mapping $x \mapsto
X_t^x(\omega)$ is continuous, differentiable, or continuously
differentiable on $\R^d$ for every $\omega \in \Omega \setminus N$
and each $t \ge 0$ given and fixed. The first spatial derivative of
the stochastic-flow coordinate $X^j$ with respect to $x_i$ will be
denoted by $\partial_i X_t^{j,x} := \partial_{x_i} X_t^{j,x} =
\partial X_t^{j,x}\! /\partial x_i$ for $t \ge 0$ and $x = (x_1,
\ldots, x_d) \in \R^d$ with $1 \le i,j \le d$. (The same notation
will also be applied to deterministic functions throughout including
their time derivatives whenever convenient.) Thus when the
stochastic flow is continuously differentiable we know that $x
\mapsto
\partial_i X_t^{j,x}(\omega)$ is continuous on $\R^d$ for every $\omega
\in \Omega \setminus N$ and each $t \ge 0$ where $P(N)=0$ and $1 \le
i,j \le d$. We will also assume that the (timewise) sample path
regularity of $X^j$ translates to the same sample path regularity of
$\partial_i X^j$, i.e.\ if $t \mapsto X_t^{j,x}(\omega)$ is
continuous or right-continuous with left limits, then $t \mapsto
\partial_i X_t^{j,x}(\omega)$ is continuous or right-continuous with
left limits for every $\omega \in \Omega \setminus N$ and each $x
\in \R^d$ where $P(N)=0$ and $1 \le i,j \le d$.

To obtain \emph{sufficient conditions} for stochastic flow
regularity, which are directly relevant for the main results in
Sections 4 and 5 below, recall that a stochastic flow $X =
(X_t^x)_{t \ge 0, x \in \R^d}$ may be viewed as a stochastic field
$Z = (Z_z)_{z \in \R_+ \times \R^d}$, where we set $Z_z = X_t^x$ for
$z = (t,x) \in \R_+ \times \R^d$, so that the results on sample path
regularity of stochastic fields are applicable to stochastic flows.
The earliest results of this kind for the existence of (H\"older)
\emph{continuous} modifications of stochastic processes (when the
index set of a stochastic field is $\R_+$) were derived by
Kolmogorov in 1934 (unpublished) and published subsequently by
Slutsky \cite{Sl} (\,see also \cite[pp.\ 158-165]{Kh} for extensions
of these results to stochastic fields when the index set is $\R_+^n$
for $n \ge 1$). Sufficient conditions for the existence of
\emph{right-continuous} modifications of stochastic processes (with
left limits) have been derived by Chentsov \cite{Ch} and Cram\'er
\cite{Cr}. Sufficient conditions for the existence of
\emph{continuously differentiable} modifications of stochastic
processes have been derived in the book by Cram\'er and Leadbetter
\cite[pp.\ 67-70]{CL}. All these conditions are of a
H\"older-in-mean type involving either two-dimensional (for
continuity) or three-dimensional (for right-continuity or
differentiability) marginal laws of the process. Different
sufficient conditions for the existence of continuously
differentiable modifications of stochastic fields (indexed by
$\R_+^n$ for $n \ge 1$) have been derived by Potthoff \cite[Theorem
3.2]{Pot} based on the ideas of Lo\`eve cited therein. These
conditions require the existence of the first partial derivative of
the original stochastic flow in the mean-square sense (thus again
being of a H\"older-in-mean type however without specifying the
admissible rate of convergence) combined with the existence of a
continuous modification of the resulting partial derivative flow
(which can be established at least formally using the extended
Kolmogorov conditions for stochastic fields referred to above).

The preceding results give a variety of general sufficient
conditions for the existence of a \emph{regular} stochastic field
and hence a \emph{regular} stochastic flow as well. Entering into a
more specific class of stochastic processes, it is well known that
SDEs driven by semimartingales with differentiable coefficients
having locally Lipschitz first partial derivatives generate
continuously differentiable flows (cf.\ \cite[Theorem 39, p.\
305]{Pr}). In particular, this is true for SDEs driven by a standard
Wiener process or a more general L\'evy process in $\R^d$. Each of
these processes therefore satisfies the hypothesis on the existence
of a continuously differentiable flow. To express the hypothesis in
a compact form we will simply say that the process $X$ can be
realised as a continuously differentiable stochastic flow
$(X_t^x)_{t \ge 0, x \in \R^d}$ in the space variable.

\vspace{6pt}

5.\ \emph{Infinitesimal generator regularity}. We will assume in the
sequel that the infinitesimal generator of $X$ is given by
\begin{align} \h{2pc} \label{2.14}
\LL_X F (x) &= \frac{1}{2} \sum_{i,j=1}^d \sigma_{ij}(x) \frac{\partial^2
F}{\partial x_i \partial x_j}(x) + \sum_{i=1}^d \mu_i(x) \frac{\partial
F}{\partial x_i} (x) - \lambda(x) F(x) \\ \notag &\h{13pt}+ \int_{\R^d
\setminus \{ 0 \}} \Big( F(y) - F(x) - \sum_{i=1}^d (y_i \m x_i) \frac{
\partial F}{\partial x_i}(x) \Big)\; \nu(x,dy)
\end{align}
for any function $F : \R^d \rightarrow \R$ from its domain and $x
\in \R^d$, where the matrix $(\sigma_{ij})_{i,j=1}^d$ with values in
$\R^{d \times d}$ is symmetric and positive semi-definite (diffusion
coefficient), the vector $(\mu_i)_{i=1}^d$ takes values in $\R^d$
(drift coefficient), $\lambda$ takes values in $\R_+$ (killing
coefficient), and $\nu(x,dy)$ is a non-negative measure on $\R^d\!
\setminus\! \{ 0 \}$ (\;\!the compensator of the measure of jumps of
$X$). For more details we refer to \cite[pp.\ 281-299]{RY} and
\cite[pp.\ 128-142]{PS}. The infinitesimal role of $\LL_X$ is
uniquely determined through its action on sufficiently regular
(smooth) functions $F$ that could also involve various boundary
conditions (\:\!on curves or surfaces in $\R^d$) depending on the
stochastic behaviour of the process $X$ (on these curves or
surfaces). It is well known that if $F$ belongs to the domain of
$\LL_X$ then
\begin{equation} \h{7pc} \label{2.15}
F(X_t) - F(X_0) - \int_0^t \LL_X F(X_s)\; ds
\end{equation}
is a (local) martingale. This is a single most useful consequence of
the previous inclusion (if known) that we will need in the sequel.
When $X$ is a semimartingale then \eqref{2.15} with $\LL_X$ from
\eqref{2.14} can also be derived for sufficiently regular (smooth)
functions $F$ using stochastic calculus techniques (It\^o's formula
and its extensions). The importance of the infinitesimal generator
\eqref{2.14} follows from the well-known fact that the optimal
stopping problem \eqref{2.1} is equivalent to the free boundary
problem
\begin{align} \h{5pc} \label{2.16}
&\LL_X V = -H\;\; \text{on}\;\; C \\[5pt] \label{2.17} &V = G
\;\; \text{on}\;\; \partial C\;\; \text{(continuous fit)} \\[3pt]
\label{2.18} &\frac{\partial V}{\partial x_i} = \frac{\partial G}
{\partial x_i}\;\; \text{on}\;\; \partial C\;\; \text{for}\;\; 1
\le i \le d\;\; \text{(smooth fit)}
\end{align}
where $\LL_X G \le -H$ on $D$, and the continuity condition
\eqref{2.17} or \eqref{2.18} applies as a variational principle when
the expectation in \eqref{2.1} with $\tau_{D'}$ in place of $\tau$
for $D' \ne D$ is discontinuous or has discontinuous first partial
derivatives at $\partial C'$ as a function of the initial point $x
\in \R^d$ respectively (for more details see \cite[p.\ 49]{PS}).
Continuity of the partial derivatives in \eqref{2.18} has been
traditionally understood/derived in the directional sense as follows
\begin{equation} \h{0pc} \label{2.19}
\lim_{h \downarrow 0} \frac{\partial V}{\partial x_i}(x_1, \ldots,
x_{i-1},x_i \pm h,x_{i+1}, \ldots, x_d) = \lim_{h \downarrow 0} \frac{
\partial G}{\partial x_i}(x_1, \ldots, x_{i-1},x_i \mp h,x_{i+1}, \ldots,
x_d)
\end{equation}
upon assuming that $(x_1, \ldots, x_{i-1},x_i \pm h,x_{i+1}, \ldots,
x_d)$ belongs to $C$ and $(x_1, \ldots, x_{i-1},x_i \mp h,x_{i+1},
\ldots, x_d)$ belongs to $D$ for $h>0$ and $1 \le i \le d$. Our main
aim in this paper is to derive the continuity of the partial
derivatives in \eqref{2.18} globally at $\partial C$, i.e.\ we aim
to show that if $x^n \in C$ converges to $x \in \partial C$ then
$(\partial V/\partial x_i)(x^n)$ converges to $(\partial G/\partial
x_i)(x)$ as $n \rightarrow \infty$ for $1 \le i \le d$. When
combined with the interior regularity results for $V$ on $C$, making
it at least continuously differentiable (in the sense of classical
derivatives), this fact will establish a global continuous
differentiability of $V$ on $\R^d$.

We will conclude this section with a few remarks on the interior
regularity of $V$ on $C$. It is well known that this can be achieved
by considering the Dirichlet/Poisson problem $\LL_X V = -H$ on a
ball (elliptic case) or a rectangle (parabolic case) contained in
$C$ where the boundary values are determined by the value function
$V$ itself upon knowing/establishing that $V$ is continuous (which
normally presents no difficulty in specific examples). Since the
boundary of a ball or a rectangle is known to be sufficiently
regular we know that the Dirichlet/Poisson problem can be solved
uniquely. For example, when $\nu \equiv 0$ in \eqref{2.14} it is
known that (locally) H\"older coefficients in \eqref{2.14} yield a
unique solution which is $C^2$ in the space variables and $C^1$ in
the time variable (see \cite[Theorem 6.13, p.\ 106]{GT} for the
elliptic case and \cite[Theorem 9, p.\ 69]{Fr-1} for the parabolic
case). This solution can then be identified with the value function
$V$ itself using the stochastic calculus or infinitesimal generator
techniques as described above (see \cite[p.\ 131]{PS} for further
details) thus establishing the interior regularity of $V$ on $C$ as
claimed. The central aim of the present paper is to establish the
$C^1$ regularity of the value function $V$ at the optimal stopping
boundary $\partial C$ that in turn is not accessible by these
arguments.

\section{Green regularity}

In this section we present two sufficient conditions for the Green
regularity of boundary points as defined in \eqref{2.8} above. The
first condition is contained in the first equivalence of
\eqref{2.11} and we expose its proof for completeness and comparison
(Lemma 1 \& Corollary 2). The second condition (Lemma 4 \& Corollary
5) has its origin in the facts that the mapping $x \mapsto \EE_x(Z)$
is finely continuous if $Z \circ \theta_t \rightarrow Z$ as $t
\downarrow 0$ where $\theta_t$ denotes the shift operator and the
implication is applicable to $Z = \sigma_U$ when $U$ is an open set
in $\R^d$ (see \cite[Corollaries 1 \& 2, p.\ 123]{Dy-1}). The two
sufficient conditions applied to stochastic flows (Corollaries 3 \&
6) will be used in the proofs of the main results in Sections 4 and
5 below.

Throughout this section we recall/assume that the (standard Markov)
process $X=(X_t)_{t \ge 0}$ and the filtration $({\cal F}_t)_{t \ge
0}$ (to which $X$ is adapted) are \emph{right-continuous} so that
the first entry and hitting times of $X$ to Borel (open and closed)
sets are stopping times (cf.\ \cite[Theorem 10.7, p.\ 54]{BG}).
Recall that $C$ denotes the continuation (open) set, $D = \R^d
\setminus C$ denotes the stopping (closed) set, and $\partial C$
denotes the boundary of the set $C$ (see Section 2 above).

\vspace{12pt}

\textbf{Lemma 1.} \emph{If $X$ is strong Feller then
\begin{equation} \h{4pc} \label{3.1}
x \mapsto \PP_{\!x}( \sigma_D \ge \eps )\;\; \text{is upper
semicontinuous on}\;\; \R^d
\end{equation}
for each $\eps>0$ given and fixed.}

\vspace{12pt}

\textbf{Proof.} Using that $\delta \p \sigma_D \circ \theta_\delta
\downarrow \sigma_D$ as $\delta \downarrow 0$, and letting $\eps>0$
be given and fixed, we find by the strong Markov property of $X$
that
\begin{align} \h{0.5pc} \label{3.2}
\PP_{\!x}(\sigma_D \ge \eps) &= \lim_{\delta \downarrow 0}\;\!
\PP_{\!x}(\delta \p \sigma_D \circ \theta_\delta \ge \eps) =
\lim_{\delta \downarrow 0}\;\! \EE_x \big[\;\! \EE_x( I(\sigma_D
\circ \theta_\delta \ge \eps \m \delta)\, \vert\, {\cal F}_\delta
)\;\! \big] \\ \notag &= \lim_{\delta \downarrow 0}\;\! \EE_x
\big[\;\! \EE_{X_\delta}( I(\sigma_D \ge \eps \m \delta) )\;\!
\big] = \lim_{\delta \downarrow 0}\;\! \EE_x \big[\;\! \PP_{\!
X_\delta}(\sigma_D \ge \eps \m \delta )\;\! \big] \\ \notag
&= \lim_{\delta \downarrow 0}\;\! \EE_x \big[\:\! F_\delta(
X_\delta) \:\!\big] = \lim_{\delta \downarrow 0}\;\! G_\delta(x)
\end{align}
where $x \mapsto F_\delta(x) := \PP_{\!x}( \sigma_D \ge \eps \m
\delta )$ is measurable so that $x \mapsto G_\delta(x) := \EE_x
\big[\:\! F_\delta(X_\delta) \:\!\big]$ is continuous on $\R^d$ by
the strong Feller property of $X$. Since moreover $\delta \mapsto
G_\delta$ is decreasing on $(0,\infty)$ as $\delta \downarrow 0$, we
see from \eqref{3.2} that \eqref{3.1} is satisfied as claimed.
\hfill $\square$

\vspace{12pt}

\textbf{Corollary 2.} \emph{If $x \in \partial C$ is
probabilistically regular for $D$ and $X$ is strong Feller, then $x$
is Green regular for $D$.}

\vspace{12pt}

\textbf{Proof.} Take any $x_n \in C$ converging to $x \in
\partial C$ as $n \rightarrow \infty$. Then by \eqref{3.1} we get
\begin{equation} \h{-1pc} \label{3.3}
0 \le \liminf_{n \rightarrow \infty}\;\! \PP_{\!x_n}(\tau_D \ge \eps)
\le \limsup_{n \rightarrow \infty}\;\! \PP_{\!x_n}(\tau_D \ge \eps)
\le \limsup_{n \rightarrow \infty}\;\! \PP_{\!x_n}(\sigma_D \ge \eps)
\le \PP_{\!x}(\sigma_D \ge \eps) = 0
\end{equation}
for each $\eps>0$ given and fixed, where the final equality follows
by probabilistic regularity of $x$ for $D$. This shows that
\eqref{2.8} is satisfied as claimed. \hfill $\square$

\vspace{12pt}

When the process $X$ can be realised as a stochastic flow
$(X_t^x)_{t \ge 0, x \in \R^d}$ we write
\begin{equation} \h{2.5pc} \label{3.4}
\tau_D^x = \inf\, \{\, t \ge 0\; \vert\; X_t^x \in D\, \}\;\;\; \&
\;\;\; \sigma_D^x = \inf\, \{\, t > 0\; \vert\; X_t^x \in D\, \}
\end{equation}
to denote the dependence of $\tau_D$ and $\sigma_D$ on $x \in \R^d$.
In this case we can reformulate the result of Corollary 2 as
follows.

\vspace{12pt}

\textbf{Corollary 3.} \emph{If $x \in \partial C$ is
probabilistically regular for $D$ and $X$ is strong Feller, then
$\tau_D^{x_n} \rightarrow 0$ in probability whenever $x_n \in C$
converges to $x \in \partial C$ as $n \rightarrow \infty$.}

\vspace{12pt}

\textbf{Proof.} This is a direct consequence of the Green regularity
established in Corollary 2. \hfill $\square$

\vspace{12pt}

When the process $X$ fails to be strong Feller then the conclusions
of Lemma 1, Corollary 2 and Corollary 3 can generally fail under
probabilistic regularity of a point from $\partial C$ for the set
$D$. We now show that the conclusions remain valid if $X$ can be
realised as a stochastic flow that is continuous in the \emph{space}
variable and a point from $\partial C$ is probabilistically regular
for the \emph{interior} $D^\circ$ of the set $D$.

\vspace{12pt}

\textbf{Lemma 4.} \emph{If $X$ can be realised as a stochastic flow
such that
\begin{equation} \h{4pc} \label{3.5}
x \mapsto X_t^x\;\; \text{is continuous on}\;\; \R^d
\end{equation}
almost surely for each $t \ge 0$ given and fixed, then
\begin{equation} \h{4pc} \label{3.6}
x \mapsto \PP_{\!x}( \sigma_{D^\circ} \ge \eps )\;\; \text{is upper
semicontinuous on}\;\; \R^d
\end{equation}
for each $\eps>0$ given and fixed.}

\vspace{12pt}

\textbf{Proof.} We first show that
\begin{equation} \h{4pc} \label{3.7}
x \mapsto \sigma_{D^\circ}^x \;\; \text{is upper
semicontinuous on}\;\; \R^d
\end{equation}
almost surely. For this, take any $x_n \rightarrow x$ in $\R^d$ as
$n \rightarrow \infty$. Denoting the exceptional set of
$\PP$\!-measure zero in \eqref{3.5} by $N_t$, and setting $N :=
\cup_{t \in \mathbb{Q}_+} N_t$ which also is a set of
$\PP$\!-measure zero, we know that \eqref{3.5} holds on $\Omega\!
\setminus\! N$ for every $t \in \mathbb{Q}_+$. Let $\omega \in
\Omega\! \setminus\! N$ be given and fixed. By definition of
$\sigma_{D^\circ}(\omega)$ and right-continuity of $t \mapsto
X_t^x(\omega)$ we know that for $\eps>0$ given and fixed, there
exists $t_\eps \in (\sigma_{D^\circ}^x(\omega), \sigma_{D^\circ}^x
(\omega) \p \eps) \cap \mathbb{Q}_+$ such that $X_{t_\eps}^x(\omega)
\in D^\circ$. Because $D^\circ$ is open it follows that there exists
$\delta_\eps>0$ such that $b(X_{t_\eps}^x(\omega),\delta_\eps)
\subseteq D^\circ$. Since \eqref{3.5} holds on $\Omega \setminus N$
for $t_\eps \in \mathbb{Q}_+$ we see that there exists $n_\eps \ge
1$ such that $X_{t_\eps}^{x_n}(\omega) \in
b(X_{t_\eps}^x(\omega),\delta_\eps)$ for all $n \ge n_\eps$. This
shows that $\sigma_{D^\circ}^{x_n}(\omega) \le t_\eps$ for all $n
\ge n_\eps$ and hence we find that $\limsup_{\,n \rightarrow \infty}
\sigma_{D^\circ}^{x_n}(\omega) \le t_\eps$. Letting $\eps \downarrow
0$ we get $\limsup_{\,n \rightarrow \infty}
\sigma_{D^\circ}^{x_n}(\omega) \le \sigma_{D^\circ}^x(\omega)$ and
this establishes \eqref{3.7} as claimed.

We next show that \eqref{3.6} holds. For this, take any $x_n
\rightarrow x$ in $\R^d$ as $n \rightarrow \infty$ and set $A_n =
\{\, \sigma_{D^\circ}^{x_n} \ge \eps\, \}$ for $n \ge 1$. Then by
Fatou's lemma for sets we find that
\begin{align} \h{1pc} \label{3.8}
\limsup_{n \rightarrow \infty}\;\! \PP_{\!x_n}(\sigma_{D^\circ} \ge \eps)
&= \limsup_{n \rightarrow \infty}\;\! \PP(\sigma_{D^\circ}^{x_n} \ge \eps)
= \limsup_{n \rightarrow \infty}\;\! \PP(A_n) \le \PP \big( \limsup_{n
\rightarrow \infty} A_n \big) \\ \notag &\le \PP(\sigma_{D^\circ}^x
\ge \eps) = \PP_{\!x} (\sigma_{D^\circ} \ge \eps)
\end{align}
where the second inequality follows since $\omega \in \limsup_{\,n
\rightarrow \infty} A_n$ if and only if $\omega \in A_{n_k}$ for $k
\ge 1$, so that $\sigma_{D^\circ}^{x_{n_k}}(\omega) \ge \eps$ for $k
\ge 1$ and hence by \eqref{3.7} we get $\sigma_{D^\circ}^x(\omega)
\ge \limsup_{\,n \rightarrow \infty} \sigma_{D^\circ}^{x_n}(\omega)
\ge \limsup_{\,k \rightarrow \infty} \sigma_{D^\circ}^{x_{n_k}}
(\omega) \ge \eps$ implying the claim. This shows that \eqref{3.6}
is satisfied as claimed. \hfill $\square$

\vspace{12pt}

\textbf{Corollary 5.} \emph{If $x \in \partial C$ is
probabilistically regular for $D^\circ$ and $X$ can be realised as a
stochastic flow such that \eqref{3.5} holds, then $x$ is Green
regular for $D^\circ$ (and thus $D$ too).}

\vspace{12pt}

\textbf{Proof.} Take any $x_n \in C$ converging to $x \in
\partial C$ as $n \rightarrow \infty$. Then similarly to the proof
of Corollary 2, we find by \eqref{3.6} that \eqref{3.3} holds with
$D^\circ$ in place of $D$ for each $\eps>0$ given and fixed, where
the final equality follows by probabilistic regularity of $x$ for
$D^\circ$. This shows that \eqref{2.8} is satisfied with $D^\circ$
in place of $D$ as claimed. \hfill $\square$

\vspace{10pt}

\textbf{Corollary 6.} \emph{If $x \in \partial C$ is
probabilistically regular for $D^\circ$ and $X$ can be realised as a
stochastic flow such that \eqref{3.5} holds, then
$\tau_{D^\circ}^{x_n} \rightarrow 0$ almost surely (and thus
$\tau_D^{x_n} \rightarrow 0$ almost surely too) whenever $x_n \in C$
converges to $x \in \partial C$ as $n \rightarrow \infty$.}

\vspace{10pt}

\textbf{Proof.} This is a direct consequence of \eqref{3.7} upon
noting that $\tau_{D^\circ}^{x_n} = \sigma_{D^\circ}^{x_n}$ for $n
\ge 1$ and $\tau_{D^\circ}^x = \sigma_{D^\circ}^x$ since $D^\circ$
is open. \hfill $\square$

\vspace{10pt}

According to \cite{Dem} and the references therein, a point $z \in
\partial C$ that is regular for $D^\circ$ is called a \emph{stable}
boundary point, and the boundary $\partial C$ is said to be
(\emph{strongly}) \emph{transversal} if $\sigma_D =
\sigma_{D^\circ}$ almost surely with respect to $\PP_x$ for all $x
\in C$ (\,for all $x \in \R^d$). Note that the results of Corollary
5 and Corollary 6 can be rephrased in terms of stable boundary
points. An important example of the strongly transversal boundary is
obtained as follows.

\vspace{10pt}

\textbf{Example 7.} If $t \mapsto b(t)$ is (piecewise)
\emph{monotone} and (left/right) \emph{continuous} on $\R_+$ and
\begin{equation} \h{6pc} \label{3.9}
D = \{\, (t,x) \in \R_+\! \times\! \R\; \vert\; x \ge b(t)\, \}
\end{equation}
then for any regular (recurrent) It\^o-McKean diffusion $X$ we have
$\sigma_D = \sigma_{D^\circ}$ almost surely with respect to $\PP_x$
for every $x=b(t)$ with $t \ge 0$ (see the proof of Corollary 8 in
\cite{CP}). Note that Corollary 5 in this case implies that
probabilistic regularity of a boundary point $z=(t,x) \in
\partial C$ implies its Green regularity despite the fact that the
time-space process $((t,X_t))_{t \ge 0}$ is not strong Feller so
that the (general) first equivalence in \eqref{2.11} is not
applicable.

\section{Continuity of the space derivative}

In this section we show that probabilistic regularity of the optimal
stopping boundary implies continuous \emph{spatial}
differentiability of the value function at the optimal stopping
boundary whenever the process admits a continuously differentiable
flow.

\vspace{4pt}

1.\ We first consider the case of infinite horizon in Theorem 8.
This will be then extended to the case of finite horizon in Theorem
10 below. Similarly to \eqref{3.4} above we write $\Lambda_t^x =
\int_0^t \lambda(X_s^x)\, ds$ to denote the dependence of
$\Lambda_t$ on $x \in \R^d$ for $t \ge 0$. We set $\delta_{i,j}=1$
if $i=j$ and $\delta_{i,j}=0$ for $i \ne j$ with $1 \le i,j \le d$.

\vspace{10pt}

\textbf{Theorem 8.} \emph{Consider the optimal stopping problem
\eqref{2.1} upon assuming that it is well posed in the sense that
the stopping time $\tau_D$ from \eqref{2.4} is optimal. Assume that
\begin{align} \h{0pc} \label{4.1}
&V\;\; \text{is continuous on}\;\; \R^d\;\; \text{and continuously
differentiable on}\;\; C\, ;\\ \label{4.2} &G\;\; \text{is continuously
differentiable on}\;\; \R^d\, ; \\ \label{4.3} &H\;\; \text{and}\;\;
\lambda\;\; \text{are Lipschitz continuous on}\;\; \R^d\;\; \text{in
the sense that} \\[6pt] \notag &\h{1.4pc} \vert H(x) \m H(y) \vert \le
K \vert x \m y \vert\;\; \&\;\; \vert \lambda(x) \m \lambda(y) \vert
\le K \vert x \m y \vert\\[6pt] \notag &\text{for all}\;\; x,y \in
\R^d \;\; \text{with some constant}\;\; K>0\;\; \text{large enough}.
\end{align}
Assume moreover that the process $X$ can be realised as a
continuously differentiable stochastic flow $(X_t^x)_{t \ge 0, x \in
\R^d}$ in the space variable and that for $z \in \partial C$ given
and fixed we have
\begin{align} \h{1pc} \label{4.4}
&\EE \bigg[\, \sup_{\alpha,\beta,\xi \in b(z,r)}\!\! e^{-\Lambda
_{\tau_D^\alpha}^\beta}\;\! \big \vert \partial_j G(X_{\tau_D
^\alpha}^\xi)\:\! \partial_i X_{\tau_D^\alpha}^{j,\xi} \big \vert
\, \bigg] < \infty \\ \label{4.5} &\EE \bigg[\, \sup_{\alpha
\in b(z,r)} \int_0^{\tau_D^\alpha}\!\! \sup_{\beta,\eta \in b
(z,r)}\!\! e^{-\Lambda_t^\beta} \big \vert \partial_i X_t^{j,
\eta} \big \vert\, dt\, \bigg] < \infty \\ \label{4.6} &\EE
\bigg[ \sup_{\alpha,\beta,\gamma \in b(z,r)}\! \bigg( e^{-\Lambda
_{\tau_D^\alpha}^\beta} \big \vert G(X_{\tau_D^\alpha}^\gamma)
\big \vert \int_0^{\tau_D^\alpha}\!\! \sup_{\eta \in b(z,r)} \big
\vert \partial_i X_t^{j,\eta} \big \vert\, dt \bigg)\;\! \bigg]
< \infty \\ \label{4.7} &\EE \bigg[\, \sup_{\alpha \in b(z,r)}
\int_0^{\tau_D^\alpha}\!\! \bigg( \sup_{\beta,\gamma \in b(z,r)}
\!\! e^{-\Lambda_t^\beta} \big \vert H(X_t^\gamma) \big \vert
\int_0^t \sup_{\eta \in b(z,r)} \big \vert \partial_i X_s^{j,\eta}
\big \vert\, ds \bigg)\, dt\, \bigg] < \infty
\end{align}
for some $r>0$ with $\partial_i X_{0+}^{j,z} = \delta_{i,j}$ for $1
\le i,j \le d$. If $X$ is strong Feller and $z$ is probabili-
stically regular for $D$, or $X$ is strong Markov and $z$ is
probabilistically regular for $D^\circ$, then
\begin{equation} \h{5pc} \label{4.8}
V\;\; \text{is continuously differentiable at}\;\; z
\end{equation}
with $\partial_i V(z) = \partial_i G(z)$ for $1 \le i \le d$. If the
hypotheses stated above hold at every $z \in \partial C$ then $V$ is
continuously differentiable on $\R^d$.}

\vspace{12pt}

\textbf{Proof.} It will be clear from the proof below that the same
arguments are applicable in any dimension $d \ge 1$ so that for ease
of notation we will assume that $d=1$ in the sequel.

\vspace{6pt}

(I): To illustrate the arguments in a clearer manner we first
consider the special case when $\Lambda_t = 0$ for $t \ge 0$. Note
that the conditions \eqref{4.6} and \eqref{4.7} are not needed in
that case.

\vspace{6pt}

1.\ Take any $x_n \in C$ converging to $z \in \partial C$ as $n
\rightarrow \infty$. Passing to a subsequence of $(x_n)_{n \ge 1}$
if needed there is no loss of generality in assuming that
\begin{equation} \h{5pc} \label{4.9}
\liminf_{n \rightarrow \infty}\, V_x(x_n) = \lim_{n \rightarrow \infty}
\frac{V(x_n \p \eps_n) \m V(x_n)}{\eps_n}
\end{equation}
for some $\eps_n \downarrow 0$ as $n \rightarrow \infty$ (we write
$V_x$ to denote $\partial V\!/\partial x$ throughout). Let $\tau_n
:= \tau_D^{x_n}$ be the optimal stopping time for $V(x_n)$ when $n
\ge 1$. Then by the mean value theorem and \eqref{4.3} we find that
\begin{align} \h{-0.5pc} \label{4.10}
V(x_n \p \eps_n) \m V(x_n) &\ge \EE \Big[ G(X_{\tau_n}^{x_n + \eps_n})
+\! \int_0^{\tau_n}\!\! H(X_t^{x_n + \eps_n})\, dt \Big]\! - \EE \Big[
G(X_{\tau_n}^{x_n}) +\! \int_0^{\tau_n}\!\! H(X_t^{x_n})\, dt \Big] \\
\notag &\h{-5pc}= \EE \Big[ G(X_{\tau_n}^{x_n + \eps_n})\! - G(X_{\tau_n}
^{x_n}) \Big]\! + \EE \Big[ \int_0^{\tau_n}\!\! \big( H(X_t^{x_n + \eps_n})
\m H(X_t^{x_n}) \big)\, dt \Big] \\ \notag &\h{-5pc}\ge \EE \Big[ G_x(X_
{\tau_n}^{\xi_n})\, \partial_x X_{\tau_n}^{\xi_n}\, \eps_n \Big]\! -
\EE \Big[ \int_0^{\tau_n}\!\! K\, \vert \partial_x X_t^{\eta_n(t)}
\vert\, \eps_n\, dt \Big]
\end{align}
where $\xi_n$ and $\eta_n(t)$ belong to $(x_n,x_n \p \eps_n)$ for $n
\ge 1$. Dividing both sides by $\eps_n$ and letting $n \rightarrow
\infty$ we find from \eqref{4.9}+\eqref{4.10} that
\begin{equation} \h{0.5pc} \label{4.11}
\liminf_{n \rightarrow \infty}\, V_x(x_n) \ge \lim_{n \rightarrow \infty}
\EE \Big[ G_x(X_{\tau_n}^{\xi_n})\, \partial_x X_{\tau_n}^{\xi_n} \Big]
- K \lim_{n \rightarrow \infty} \EE \Big[ \int_0^{\tau_n}\!\! \vert
\partial_x X_t^{\eta_n(t)} \vert\, dt \Big] = G_x(z)
\end{equation}
where in the final equality we use that $\tau_n \rightarrow 0$
almost surely as $n \rightarrow \infty$ by Green regularity of $z$
for $D$ as established in Corollary 3 and Corollary 6 above (in the
former case one may need to pass to a subsequence of $(x_n)_{n \ge
1}$ which is sufficient for the present purposes) combined with the
dominated convergence theorem which is applicable due to \eqref{4.4}
and \eqref{4.5} respectively.

\vspace{6pt}

2.\ Similarly, there is no loss of generality in assuming that
\begin{equation} \h{5pc} \label{4.12}
\limsup_{n \rightarrow \infty}\, V_x(x_n) = \lim_{n \rightarrow \infty}
\frac{V(x_n) \m V(x_n \m \eps_n)}{\eps_n}
\end{equation}
for some $\eps_n \downarrow 0$ as $n \rightarrow \infty$. By the
mean value theorem and \eqref{4.3} we find that
\begin{align} \h{-0.5pc} \label{4.13}
V(x_n) \m V(x_n \m \eps_n) &\le \EE \Big[ G(X_{\tau_n}^{x_n}) +\! \int_0^
{\tau_n}\!\! H(X_t^{x_n})\, dt \Big]\! - \EE \Big[ G(X_{\tau_n}^{x_n -
\eps_n}) +\! \int_0^{\tau_n}\!\! H(X_t^{x_n - \eps_n})\, dt \Big] \\
\notag &= \EE \Big[ G(X_{\tau_n}^{x_n})\! - G(X_{\tau_n}^{x_n - \eps_n})
\Big]\! + \EE \Big[ \int_0^{\tau_n}\!\! \big( H(X_t^{x_n}) \m H(X_t^
{x_n - \eps_n}) \big)\, dt \Big] \\ \notag &\le \EE \Big[ G_x(X_
{\tau_n}^{\xi_n})\, \partial_x X_{\tau_n}^{\xi_n}\, \eps_n \Big]\!
+ \EE \Big[ \int_0^{\tau_n}\!\! K\, \vert \partial_x X_t^{\eta_n(t)}
\vert\, \eps_n\, dt \Big]
\end{align}
where $\xi_n$ and $\eta_n(t)$ belong to $(x_n \m \eps_n,x_n)$ for $n
\ge 1$. Dividing both sides by $\eps_n$ and letting $n \rightarrow
\infty$ we find from \eqref{4.12}+\eqref{4.13} that
\begin{equation} \h{0.5pc} \label{4.14}
\limsup_{n \rightarrow \infty}\, V_x(x_n) \le \lim_{n \rightarrow \infty}
\EE \Big[ G_x(X_{\tau_n}^{\xi_n})\, \partial_x X_{\tau_n}^{\xi_n} \Big]
+ K \lim_{n \rightarrow \infty} \EE \Big[ \int_0^{\tau_n}\!\! \vert
\partial_x X_t^{\eta_n(t)} \vert\, dt \Big] = G_x(z)
\end{equation}
where in the final equality we use the same arguments as following
\eqref{4.11} above. Combining \eqref{4.11} and \eqref{4.14} we see
that $\lim_{\,n \rightarrow \infty} V_x(x_n) = G_x(z)$ and this
completes the proof when $\Lambda_t = 0$ for $t \ge 0$.

\vspace{6pt}

(II): Next we consider the general case when $\Lambda_t \ne 0$ for
$t \ge 0$. Note that the conditions \eqref{4.6} and \eqref{4.7} are
needed in that case unless $\lambda$ is constant for all $t \ge 0$.
The proof in the general case can be carried out along the same
lines as in the special case above and we only highlight the needed
modifications throughout.

\vspace{6pt}

3.\ Taking any $x_n \in C$ converging to $z \in \partial C$ as $n
\rightarrow \infty$ and arguing as in \eqref{4.9} above, we see that
the right-hand side of the first inequality in \eqref{4.10} reads as
follows
\begin{align} \h{1pc} \label{4.15}
&\EE \Big[ e^{-\Lambda_{\tau_n}^{x_n+\eps_n}} G(X_{\tau_n}^{x_n+
\eps_n}) - e^{-\Lambda_{\tau_n}^{x_n}} G(X_{\tau_n}^{x_n}) \Big]
\\ \notag &\h{13pt}+ \EE \Big[ \int_0^{\tau_n} \Big( e^{-\Lambda
_t^{x_n+\eps_n}} H(X_t^{x_n+\eps_n}) - e^{-\Lambda_t^{x_n}} H(X
_t^{x_n}) \Big)\, dt\, \Big] \\[2pt] \notag &= \EE \Big[ e^{-
\Lambda_{\tau_n}^{x_n}} \Big( e^{\Lambda_{\tau_n}^{x_n}-\Lambda
_{\tau_n}^{x_n+\eps_n}}\!\! - 1 \Big) G(X_{\tau_n}^{x_n+\eps_n})
\Big] + \EE \Big[ e^{-\Lambda_{\tau_n}^{x_n}} \Big( G(X_{\tau_n}
^{x_n+\eps_n}) - G(X_{\tau_n}^{x_n}) \Big) \Big] \\[1pt] \notag
&\h{13pt}+ \EE \Big[ \int_0^{\tau_n} e^{-\Lambda_t^{x_n}} \Big(
e^{\Lambda_t^{x_n}-\Lambda_t^{x_n+\eps_n}}\!\! - 1 \Big) H(X_t
^{x_n+\eps_n})\, dt\, \Big] \\ \notag &\h{13pt}+ \EE \Big[
\int_0^{\tau_n} e^{-\Lambda_t^{x_n}} \Big( H(X_t^{x_n+\eps_n})
- H(X_t^{x_n}) \Big)\, dt\, \Big]
\end{align}
for $n \ge 1$. The second expectation and the fourth expectation on
the right-hand side of \eqref{4.15} can be handled in exactly the
same way as the corresponding two expectations in \eqref{4.10}, and
this yields the conclusion of \eqref{4.11} above, i.e.
\begin{equation} \h{7pc} \label{4.16}
\liminf_{n \rightarrow \infty}\, V_x(x_n) \ge G_x(z)
\end{equation}
provided that the liminf of the first expectation on the right-hand
side of \eqref{4.15} divided by $\eps_n$ and the liminf of the third
expectation on the right-hand side of \eqref{4.15} divided by
$\eps_n$ are non-negative as $n \rightarrow \infty$. To see that
both liminfs are non-negative, note that \eqref{4.3} and the mean
value theorem imply that
\begin{align} \h{3pc} \label{4.17}
\frac{e^{\Lambda_{\sigma_n}^{x_n} - \Lambda_{\sigma_n}^{x_n+\eps_n}}-1}
{\eps_n} &\ge \frac{e^{-\eps_n K \int_0^{\sigma_n} \vert \partial_x X_s
^{\eta_n(s)} \vert\, ds}-1}{\eps_n} \\ \notag &= -K \big( \textstyle
\int_0^{\sigma_n} \vert \partial_x X_s^{\eta_n(s)} \vert\, ds \big)\;\!
e^{-\zeta_n K \int_0^{\sigma_n} \vert \partial_x X_s^{\eta_n(s)}
\vert\, ds}
\end{align}
with $\sigma_n$ equal to either $\tau_n$ (the first expectation) or
$t \in [0,\tau_n]$ (the third expectation) where $\eta_n(s)$ belongs
to $(x_n,x_n \p \eps_n)$ and $\zeta_n$ belongs to $(0,\eps_n)$ for
$n \ge 1$. Using then the same arguments as in \eqref{4.11} above
with \eqref{4.6}+\eqref{4.7} in place of \eqref{4.4}+\eqref{4.5}, we
see that the inequality \eqref{4.17} yields the fact that the two
liminfs are non-negative so that \eqref{4.16} holds as claimed.

\vspace{6pt}

4.\ Similarly, arguing as in \eqref{4.12} we see that the right-hand
side of the first inequality in \eqref{4.13} reads as follows
\begin{align} \h{1pc} \label{4.18}
&\EE \Big[ e^{-\Lambda_{\tau_n}^{x_n}} G(X_{\tau_n}^{x_n}) -e^{
-\Lambda_{\tau_n}^{x_n-\eps_n}} G(X_{\tau_n}^{x_n-\eps_n}) \Big]
\\ \notag &\h{13pt}+ \EE \Big[ \int_0^{\tau_n} \Big( e^{-\Lambda
_t^{x_n}} H(X_t^{x_n}) - e^{-\Lambda_t^{x_n-\eps_n}} H(X_t^{x_n
-\eps_n}) \Big)\, dt\, \Big] \\[2pt] \notag &= \EE \Big[ e^{-
\Lambda_{\tau_n}^{x_n}} \Big( 1 - e^{\Lambda_{\tau_n}^{x_n}-
\Lambda_{\tau_n}^{x_n-\eps_n}} \Big) G(X_{\tau_n}^{x_n}) \Big]
+ \EE \Big[ e^{-\Lambda_{\tau_n}^{x_n-\eps_n}} \Big( G(X_{\tau_n}
^{x_n}) - G(X_{\tau_n}^{x_n-\eps_n}) \Big) \Big] \\[1pt] \notag
&\h{13pt}+ \EE \Big[ \int_0^{\tau_n} e^{-\Lambda_t^{x_n}} \Big(
1 - e^{\Lambda_t^{x_n}-\Lambda_t^{x_n-\eps_n}} \Big) H(X_t
^{x_n})\, dt\, \Big] \\ \notag &\h{13pt}+ \EE \Big[
\int_0^{\tau_n} e^{-\Lambda_t^{x_n-\eps_n}} \Big( H(X_t^{x_n})
- H(X_t^{x_n-\eps_n}) \Big)\, dt\, \Big]
\end{align}
for $n \ge 1$. The second expectation and the fourth expectation on
the right-hand side of \eqref{4.18} can be handled in exactly the
same way as the corresponding two expectations in \eqref{4.13} and
this yields the conclusion of \eqref{4.14} above, i.e.
\begin{equation} \h{7pc} \label{4.19}
\limsup_{n \rightarrow \infty}\, V_x(x_n) \le G_x(z)
\end{equation}
provided that the limsup of the first expectation on the right-hand
side of \eqref{4.18} divided by $\eps_n$ and the limsup of the third
expectation on the right-hand side of \eqref{4.18} divided by
$\eps_n$ are non-positive as $n \rightarrow \infty$. To see that
both limsups are non-positive, note that \eqref{4.3} and the mean
value theorem imply that
\begin{align} \h{3pc} \label{4.20}
\frac{1 - e^{\Lambda_{\sigma_n}^{x_n}-\Lambda_{\sigma_n}^{x_n-\eps_n}}}
{\eps_n} &\le \frac{1 - e^{-\eps_n K \int_0^{\sigma_n} \vert \partial_x X_s
^{\eta_n(s)} \vert\, ds}}{\eps_n} \\ \notag &\h{-5pc}= K \big( \textstyle
\int_0^{\sigma_n} \vert \partial_x X_s^{\eta_n(s)} \vert\, ds \big)\;\!
e^{-\zeta_n K \int_0^{\sigma_n} \vert \partial_x X_s^{\eta_n(s)}
\vert\, ds}
\end{align}
with $\sigma_n$ equal to either $\tau_n$ (the first expectation) or
$t \in [0,\tau_n]$ (the third expectation) where $\eta_n(s)$ belongs
to $(x_n \m \eps_n,x_n)$ and $\zeta_n$ belongs to $(0,\eps_n)$ for
$n \ge 1$. Using then the same arguments as in \eqref{4.14} above
with \eqref{4.6}+\eqref{4.7} in place of \eqref{4.4}+\eqref{4.5}, we
see that the inequality \eqref{4.20} yields the fact that the two
limsups are non-positive so that \eqref{4.19} holds as claimed.
Combining \eqref{4.16} and \eqref{4.19} we see that $\lim_{\,n
\rightarrow \infty} V_x(x_n) = G_x(z)$ and this completes the proof
when $\Lambda_t \ne 0$ for $t \ge 0$. \hfill $\square$

\vspace{16pt}

\textbf{Remark 9.} Note that the conditions \eqref{4.4}-\eqref{4.7}
are used in the proof above as sufficient conditions for the
dominated convergence theorem to establish the convergence relations
\eqref{4.11} and \eqref{4.14} (when $\lambda$ is zero) and their
extensions \eqref{4.16} and \eqref{4.19} (when $\lambda$ is not
constant). (Recall from the proof that the conditions \eqref{4.6}
and \eqref{4.7} are not needed when $\lambda$ is constant.) These
sufficient conditions, although applicable in a large number of
examples, are not necessary in general and in some specific examples
one can often exploit additional information (e.g.\ the
geometric/analytic structure of the optimal stopping boundary) and
derive the convergence relations without appealing to the dominated
convergence theorem (see the proof of Theorem 3.1 in \cite{Pe-1} for
such an example). As it is exceedingly complicated to describe all
possible ways that lead to relaxed forms of the sufficient
conditions \eqref{4.4}-\eqref{4.7}, we have stated them in their
present form with a view that the structure of the proof above
remains unchanged if these sufficient conditions are replaced by
other/weaker ones. A similar remark applies to the condition
\eqref{4.3}. For instance, replacing the global Lipschitz continuity
of $H$ in \eqref{4.3} by a local Lipschitz continuity in the sense
that
\begin{equation} \h{7pc} \label{4.21}
\vert H(x) \m H(y) \vert \le K_n \vert x \m y \vert
\end{equation}
for all $x,y \in b(z,R_n)$ with some constant $K_n>0$ large enough
where $R_n \rightarrow \infty$ as $n \rightarrow \infty$, it is seen
from the proof above that the result of Theorem 8 (\:\!with
$\lambda=0$) remains valid if
\begin{equation} \h{7pc} \label{4.22}
\lim_{n \rightarrow \infty} K_n\, \EE \Big[ \int_0^{\bar \tau_n}\!\!
\vert \partial_x X_t^{\eta_n} \vert\, dt \Big] = 0
\end{equation}
where $\bar \tau_n := \tau_n \wedge \inf\, \{\, t \ge 0\; \vert\;
X_t^{x_n + \eps_n} \notin b(z,R_n)\;\; \text{or}\;\; X_t^{x_n}
\notin b(z,R_n)\, \}$ and $R_n>0$ is chosen large enough so that
\begin{equation} \h{1pc} \label{4.23}
\EE \Big[ \int_0^{\bar \tau_n}\!\! \big \vert H(X_t^{x_n \pm \eps_n}) \m
H(X_t^{x_n}) \big \vert\, dt \Big] \ge \EE \Big[ \int_0^{\tau_n}\!\! \big
\vert H(X_t^{x_n \pm \eps_n}) \m H(X_t^{x_n}) \big \vert\, dt \Big] -
\eps_n\:\! \delta
\end{equation}
for all $n \ge 1$ with $\delta>0$ given and fixed. Similarly, the
global Lipschitz continuity of $\lambda$ in \eqref{4.3} can be
replaced by a local Lipschitz continuity and we will omit further
details. Finally, the proof above shows that it is sufficient to
have continuous differentiability of the flow near the optimal
stopping boundary only.

\vspace{12pt}

2.\ The optimal stopping problem \eqref{2.1} considered in Theorem 8
has infinite horizon. The arguments used in the proof carry over to
the optimal stopping problem \eqref{2.2} with finite horizon as long
as continuous \emph{spatial} differentiability of the value function
is considered. We formally present this extension in the next
theorem. Continuous \emph{temporal} differentiability of the value
function requires different arguments and will be considered in the
next section.

Recall that the optimal stopping problem \eqref{2.2} includes the
case when the functions $G$ and $H$ are time dependent which can be
formally obtained by setting $X_t^1 = t$ for $t \ge 0$. Thus the
process $X$ in this case is given by $X_t=(t,X_t^2, \ldots, X_t^d)$
for $t \ge 0$. The continuation set is given by $C = \{\, (t,x) \in
[0,T]\! \times\! \R^{d-1}\; \vert\; V(t,x) > G(t,x)\, \}$ and the
stopping set is given by $D = \{\, (t,x) \in [0,T]\! \times\!
\R^{d-1}\; \vert\; V(t,x) = G(t,x)\, \}$. Note that the process $C
:= X^1$ can always be realised as a stochastic flow by setting
$C_s^t = t \p s$ for $t \ge 0$ and $s \ge 0$. Hence when $(X^2,
\ldots, X^d)$ can be realised as a stochastic flow in the space
variable $x$ from $\R^{d-1}$ we will denote the entire flow by
$(X_s^{t,x})$ for $s \ge 0$ and $(t,x) \in [0,T] \times \R^{d-1}$.
Note that $X_0^{t,x} = (t,x_2, \ldots,x_d)$ for $t \in [0,T]$ and $x
= (x_2, \ldots, x_d) \in \R^{d-1}$.

\vspace{12pt}

\textbf{Theorem 10.} \emph{Consider the optimal stopping problem
\eqref{2.2} upon assuming that it is well posed in the sense that
the stopping time $\tau_D$ from \eqref{2.4} is optimal. Assume that
\begin{align} \h{0pc} \label{4.24}
&\h{-8pt}V\;\, \text{is continuous on}\;\, [0,T]\! \times\! \R^{d-1}
\, \text{and continuously differentiable on}\;\, C\, ;\\ \label{4.25}
&\h{-8pt}G\;\, \text{is continuously differentiable on}\;\, [0,T]\!
\times\! \R^{d-1}\, ; \\ \label{4.26} &\h{-8pt}x \mapsto H(t,x)\;
\, \text{and}\;\, x \mapsto \lambda(t,x)\;\, \text{are Lipschitz
continuous on}\;\, \R^{d-1}\;\, \text{in the sense that} \\[6pt]
\notag &\h{2pc} \vert H(t,x) \m H(t,y) \vert \le K \vert x \m y
\vert\;\; \&\;\; \vert \lambda(t,x) \m \lambda(t,y) \vert \le K
\vert x \m y \vert \\[6pt] \notag &\h{-8pt}\text{for every}\;\;
t \in [0,T]\;\; \text{and all}\;\; x,y \in \R^{d-1}\;\,
\text{with some constant} \;\; K>0\;\; \text{large enough}.
\end{align}
Assume moreover that the process $X$ can be realised as a
continuously differentiable stochastic flow $(X_s^{t,x})$ in the
space variable for $s \ge 0$ and $(t,x) \in [0,T] \times \R^{d-1}$
and that for $z \in
\partial C$ given and fixed the conditions \eqref{4.4}-\eqref{4.7}
are satisfied for some $r>0$ with $\partial_i X_{0+}^{j,z} =
\delta_{i,j}$ for $2 \le i \le d$ and $1 \le j \le d$. If $z$ is
probabilistically regular for $D^\circ$ then
\begin{equation} \h{4pc} \label{4.27}
\partial_2 V, \ldots, \partial_d V\;\; \text{exist and are continuous
at}\;\; z
\end{equation}
with $\partial_i V(z) = \partial_i G(z)$ for $2 \le i \le d$. If the
hypotheses stated above hold at every $z \in \partial C$ then
$\partial_2 V, \ldots, \partial_d V$ exist and are continuous on
$[0,T]\! \times\! \R^{d-1}$.}

\vspace{12pt}

\textbf{Proof.} This can be established using exactly the same
arguments as in the proof of Theorem 8 upon noting that adding
$\eps_n$ to any but the first (time) coordinate of the process $X$
does not alter the remaining time horizon. \hfill $\square$

\vspace{16pt}

\textbf{Remark 11.} Note that the comments on the sufficient
conditions from Theorem 8 made in Remark 9 above extend to the
corresponding sufficient conditions in Theorem 10 and we will omit
further details in this direction.

\vspace{12pt}

3.\ The result and proof of Theorems 8 and 10 extend to the case
when the gain function $G$ in the optimal stopping problem
\eqref{2.1}/\eqref{2.2} is \emph{not smooth} away from the optimal
stopping boundary $\partial C$. Instead of formulating a general
theorem of this kind, which would be overly technical and rather
difficult to read, we will illustrate key arguments of such
extensions through an important example next. A different method of
proof is based on extensions of the It\^o-Tanaka formula dealing
with singularities of $G$ on curves and surfaces (cf.\ \cite{Pe-2}
and \cite{Pe-3}) and this will be presented in the next section.

\vspace{12pt}

\np

\textbf{Example 12 (Continuity of the space derivative in the
American put)}. Consider the optimal stopping problem
\begin{equation} \h{6pc} \label{4.28}
V(t,x) = \sup_{0 \le \tau \le T-t} \EE \Big[ e^{-r \tau} \big( K\!
- X_\tau^x \big)^+ \Big]
\end{equation}
where $(t,x) \in [0,T]\! \times\! (0,\infty)$, $r>0$, $K>0$ and the
supremum is taken over stopping times $\tau$ of $X$ solving the
stochastic differential equation
\begin{equation} \h{7pc} \label{4.29}
d X_t = r X_t\, dt + \sigma\:\! X_t\, dB_t
\end{equation}
with $X_0 = x$ where $\sigma>0$ and $B$ is a standard Brownian
motion (see \cite[Section 25]{PS} for further details). Horizon in
the optimal stopping problem \eqref{4.28} is finite so that the
setting belongs to Theorem 10 above. Since the gain function $G(x)
:= (K \m x)^+$ for $x>0$ is not differentiable at $K$ we see that
the condition \eqref{4.25} fails and hence we cannot conclude that
\begin{equation} \h{7pc} \label{4.30}
V_x\;\; \text{is continuous on}\;\; \partial C
\end{equation}
using Theorem 10 (we write $V_x$ to denote $\partial V\!/\partial x$
throughout). We will now show however that the method of proof of
Theorems 8 and 10 extends to cover the case of the
non-differentiable gain function $G(x) = (K \m x)^+$ for $x>0$. This
will also serve as an illustration of how similar other cases of
non-smooth gain functions $G$ in the optimal stopping problem
\eqref{2.1}/\eqref{2.2} can be handled. The derivation of
\eqref{4.30} will be divided in three steps as follows.

\vspace{6pt}

1.\ Well-known arguments show that the optimal stopping time in
\eqref{4.28} equals $\tau_D^{t,x} = \inf\, \{\, s \in [0,T \m t]\;
\vert\; X_s^x \le b(t \p s)\, \}$ where the optimal stopping
boundary $t \mapsto b(t)$ is increasing on $[0,T]$ with $0 < b(0) <
b(T) = K$ (see \cite[Subsection 25.2]{PS}). If a point $z = (t,b(t))
\in
\partial C$ is given and fixed, then by the increase of $b$ combined
with the law of iterated logarithm for standard Brownian motion
(cf.\ \cite[p.\ 112]{KS}) we see that $z$ is probabilistically
regular for $D^\circ$ (formally this could also be derived from
probabilistic regularity of $z$ for $D$ combined with the fact of
Example 7 above). Since $X$ can be realised as a continuous
stochastic flow $x \mapsto x X_t^1$ on $(0,\infty)$, where we set
$X_t^1 = \exp(\sigma B_t \p (r \m \sigma^2\!/2)\;\! t)$ for $t \ge
0$, it follows by Corollary 5 that $z$ is Green regular for
$D^\circ$. Taking any sequence $(t_n,x_n) \in C$ converging to $z$
as $n \rightarrow \infty$, it follows therefore by Corollary 6 that
$\tau_D^{t_n,x_n} \rightarrow 0$ almost surely as $n \rightarrow
\infty$. Note that the latter Green regularity has been obtained
without appeal to a strong Feller property which fails for the
time-space process $((t,X_t))_{0 \le t \le T}$ in this case.

\vspace{6pt}

2.\ We next connect to the first part of the proof of Theorems 8 and
10. Passing to a subsequence of $((t_n,x_n))_{n \ge 1}$ if needed
there is no loss of generality in assuming that
\begin{equation} \h{3.5pc} \label{4.31}
\liminf_{n \rightarrow \infty}\, V_x(t_n,x_n) = \lim_{n \rightarrow
\infty} \frac{V(t_n,x_n \p \eps_n) \m V(t_n,x_n)}{\eps_n}
\end{equation}
for some $\eps_n \downarrow 0$ as $n \rightarrow \infty$. Let
$\tau_n := \tau_D^{t_n,x_n}$ denote the optimal stopping time for
$V(t_n,x_n)$ when $n \ge 1$. Then using that $K > x_n X_{\tau_n}^1$
if and only if $\tau_n < T \m t_n$ we find that
\begin{align} \h{0pc} \label{4.32}
&V(t_n,x_n \p \eps_n) \m V(t_n,x_n) \ge \EE \Big[ e^{-r \tau_n}
\big( K \m (x_n \p \eps_n) X_{\tau_n}^1 \big)^+ \Big] - \EE \Big[
e^{-r \tau_n} \big( K \m x_n X_{\tau_n}^1 \big)^+ \Big] \\[3pt]
\notag &\h{3pc}\ge \EE \Big[ \Big( e^{-r \tau_n} \big( K \m (x_n
\p \eps_n) X_{\tau_n}^1 \big) - e^{-r \tau_n} \big( K \m x_n
X_{\tau_n}^1 \big) \Big)\:\! I(\tau_n < T \m t_n) \Big] \\[3pt]
\notag &\h{5pc}= \EE \Big[ e^{-r \tau_n} (-\eps_n)\:\! X_{\tau_n}^1
\:\! I(\tau_n < T \m t_n) \Big]
\end{align}
for $n \ge 1$. Dividing both sided by $\eps_n$ and letting $n
\rightarrow \infty$ we find from \eqref{4.31}+\eqref{4.32} that
\begin{equation} \h{3pc} \label{4.33}
\liminf_{n \rightarrow \infty}\, V_x(t_n,x_n) \ge - \lim_{n \rightarrow
\infty} \EE \Big[ e^{-r \tau_n} X_{\tau_n}^1 \:\! I(\tau_n < T \m t_n)
\Big] = -1
\end{equation}
where in the last equality we use that $\tau_n \rightarrow 0$ almost
surely as $n \rightarrow \infty$ combined with the dominated
convergence theorem due to $\EE(\sup_{\,0 \le t \le T} X_t^1) <
\infty$.

\vspace{6pt}

3.\ We finally connect to the second part of the proof of Theorems 8
and 10. Similarly, there is no loss of generality in assuming that
\begin{equation} \h{3pc} \label{4.34}
\limsup_{n \rightarrow \infty}\, V_x(t_n,x_n) = \lim_{n \rightarrow
\infty} \frac{V(t_n,x_n) \m V(t_n,x_n \m \eps_n)}{\eps_n}
\end{equation}
for some $\eps_n \downarrow 0$ as $n \rightarrow \infty$. Then using
the same arguments as in \eqref{4.32} we find that
\begin{align} \h{0pc} \label{4.35}
&V(t_n,x_n) \m V(t_n,x_n \m \eps_n) \le \EE \Big[ e^{-r \tau_n}
\big( K \m x_n X_{\tau_n}^1 \big)^+ \Big] - \EE \Big[ e^{-r \tau_n}
\big( K \m (x_n \m \eps_n) X_{\tau_n}^1 \big)^+ \Big] \\[3pt]
\notag &\h{3pc}\le \EE \Big[ \Big( e^{-r \tau_n} \big( K \m x_n
X_{\tau_n}^1 \big) - e^{-r \tau_n} \big( K \m (x_n \m \eps_n)
X_{\tau_n}^1 \big) \Big)\:\! I(\tau_n < T \m t_n) \Big] \\[3pt]
\notag &\h{3pc}= \EE \Big[ e^{-r \tau_n} (-\eps_n)\:\! X_{\tau_n}^1
\:\! I(\tau_n < T \m t_n) \Big]
\end{align}
for $n \ge 1$. Dividing both sided by $\eps_n$ and letting $n
\rightarrow \infty$ we find from \eqref{4.34}+\eqref{4.35} that
\begin{equation} \h{3pc} \label{4.36}
\limsup_{n \rightarrow \infty}\, V_x(t_n,x_n) \le - \lim_{n \rightarrow
\infty} \EE \Big[ e^{-r \tau_n} X_{\tau_n}^1 \:\! I(\tau_n < T \m t_n)
\Big] = -1
\end{equation}
where in the last equality we use the same arguments as in
\eqref{4.33} above. Combining \eqref{4.33} and \eqref{4.36} we see
that $\lim_{\,n \rightarrow \infty} V_x(t_n,x_n) = G_x(z) = -1$ and
this completes the proof of \eqref{4.30}.

\section{Continuity of the time derivative}

In this section we show that probabilistic regularity of the optimal
stopping boundary implies continuous temporal differentiability of
the value function at the optimal stopping boundary whenever the
process admits a continuous flow. We assume throughout that the
process is given by $X_t = (t,X_t^2, \ldots ,X_t^d)$ for $t \ge 0$
as discussed prior to Theorem 10 above.

\vspace{6pt}

1.\ We first consider the case of infinite horizon in Theorem 13.
This will be then extended to the case of finite horizon in Theorem
15 below.

\vspace{12pt}

\textbf{Theorem 13.} \emph{Consider the optimal stopping problem
\eqref{2.1} upon assuming that it is well posed in the sense that
the stopping time $\tau_D$ from \eqref{2.4} is optimal. Assume that
\begin{align} \h{0pc} \label{5.1}
&\h{-8pt} V\;\; \text{is continuous on}\;\; \R_+\! \times\! \R^{d-1}\;\;
\text{and continuously differentiable on}\;\; C\, ;\\ \label{5.2} &\h{-8pt}
G\;\; \text{is continuously differentiable on}\;\; \R_+\! \times\! \R^{d-1}
\, ; \\ \label{5.3} &\h{-8pt} t \mapsto H(t,x)\;\; \text{and}\;\; t \mapsto
\lambda(t,x)\;\; \text{are Lipschitz continuous on}\;\; \R_+\;\; \text{in
the sense that} \\[6pt] \notag &\h{2pc} \vert H(t,x) \m H(s,x) \vert \le
K \vert t \m s \vert\;\; \&\;\; \vert \lambda(t,x) \m \lambda(s,x) \vert
\le K \vert t \m s \vert\\[6pt] \notag &\h{-8pt}\text{for all}\;\; t,s
\in \R_+\;\; \text{and every}\;\; x \in \R^{d-1}\; \text{with some
constant} \;\; K>0\;\; \text{large enough}.
\end{align}
Assume moreover that the process $X$ can be realised as a continuous
stochastic flow $(X_s^{t,x})$ in the space variable for $s \ge 0$
and $(t,x) \in \R_+\! \times\! \R^{d-1}$ and that for $z \in
\partial C$ given and fixed the following conditions are
satisfied
\begin{align} \h{4pc} \label{5.4}
&\EE \bigg[\, \sup_{\alpha,\beta,\xi \in b(z,r)}\!\! e^{-\Lambda
_{\tau_D^\alpha}^\beta}\;\! \big \vert \partial_t G(X_{\tau_D
^\alpha}^\xi) \big \vert\, \bigg] < \infty \\ \label{5.5} &\EE
\bigg[\, \sup_{\alpha \in b(z,r)} \int_0^{\tau_D^\alpha}\!\!
\sup_{\beta \in b(z,r)}\!\! e^{-\Lambda_t^\beta}\, dt\, \bigg]
< \infty \\ \label{5.6} &\EE \bigg[\, \bigg( \sup_{\alpha,\beta,
\gamma \in b(z,r)}\!\! e^{-\Lambda_{\tau_D^\alpha}^\beta} \big
\vert G(X_{\tau_D^\alpha}^\gamma)\big \vert\, \tau_D^\alpha\:\!
\bigg)\, \bigg] < \infty \\ \label{5.7} &\EE \bigg[\, \sup_{\alpha
\in b(z,r)} \int_0^{\tau_D^\alpha}\!\! \bigg( \sup_{\beta,\gamma
\in b(z,r)}\!\! e^{-\Lambda_t^\beta} \big \vert H(X_t^\gamma)
\big \vert\, t\:\! \bigg)\, dt\, \bigg] < \infty
\end{align}
for some $r>0$. If $z$ is probabilistically regular for $D^\circ$
then
\begin{equation} \h{4pc} \label{5.8}
\partial_t V\;\; \text{exists and is continuous at}\;\; z
\end{equation}
with $\partial_t V(z) = \partial_t G(z)$. If the hypotheses stated
above hold at every $z \in \partial C$ then $\partial_t V$ exists
and is continuous on $\R_+\! \times\! \R^{d-1}$.}

\vspace{12pt}

\textbf{Proof.} Due to $X_t = (t,X_t^2, \ldots ,X_t^d)$ for $t \ge
0$ as assumed throughout we see that the setting of Theorem 13
reduces to the setting of Theorem 8. All the claims therefore follow
by applying Theorem 8 upon noting that $\partial_t X_t^{1,z} = 1$
and $\partial_t X_t^{i,z} = 0$ for $2 \le i \le d$ with $t \ge 0$
and $z \in \R_+\! \times\! \R^{d-1}$ so that the sufficient
conditions \eqref{4.4}-\eqref{4.7} in Theorem 8 transform to the
sufficient conditions \eqref{5.4}-\eqref{5.7} stated above. \hfill
$\square$

\vspace{16pt}

\textbf{Remark 14.} Note that the comments on the sufficient
conditions from Theorem 8 made in Remark 9 above extend to the
corresponding sufficient conditions in Theorem 13 and we will omit
further details in this direction.

\vspace{12pt}

2.\ The optimal stopping problem considered in Theorem 13 has
infinite horizon and the arguments used in the proof are analogous
to the arguments used in the proofs of Theorems 8 and 10 above.
Continuous temporal differentiability of the value function on
finite horizon requires different arguments and will be considered
in the next theorem. A key difficulty in the previous approach is
that adding $\eps_n$ to the first (time) coordinate of the process
$X$ (see \eqref{4.10} above) alters the remaining time horizon so
that the stopping time which is \emph{optimal} for $V(t_n,x_n)$ is
\emph{no longer} admissible for $V(t_n \p \eps_n,x_n)$ with $n \ge
1$. To overcome this difficulty we will apply a Taylor expansion of
the second order (It\^o's formula) instead of the first order as in
the proofs of Theorems 8 and 10 above.

\vspace{12pt}

\textbf{Theorem 15.} \emph{Consider the optimal stopping problem
\eqref{2.2} upon assuming that it is well posed in the sense that
the stopping time $\tau_D$ from \eqref{2.4} is optimal. Assume that
\begin{align} \h{0pc} \label{5.9}
&\h{-8pt}V\;\, \text{is continuous on}\;\, [0,T]\! \times\! \R^{d-1}
\, \text{and continuously differentiable on}\;\, C\, ; \\[3pt] \label
{5.10} &\h{-8pt}(t,x) \mapsto G(t,x)\;\, \text{is once continuously
differentiable with respect to}\;\, t\;\; \text{and twice} \\[-3pt]
\notag &\h{-8pt}\text{continuously differentiable with respect to}
\;\; x\;\; \text{on}\;\; [0,T]\! \times\! \R^{d-1}\, ; \\[3pt] \label
{5.11} &\h{-8pt}t \mapsto \tilde H(t,x) := (G_t \p \LL_X G \p H)(t,x)
\;\, \text{and}\;\; t \mapsto \lambda(t,x)\;\; \text{are Lipschitz
continuous} \\[-3pt] \notag &\h{-8pt}\text{on}\;\; [0,T]\;\; \text
{in the sense that} \\[6pt] \notag &\h{2pc} \vert \tilde H(t,x) \m
\tilde H(s,x) \vert \le K \vert t \m s \vert\;\; \&\;\; \vert \lambda
(t,x) \m \lambda(s,x) \vert \le K \vert t \m s \vert \\[6pt] \notag
&\h{-8pt} \text{for all}\;\; t,s \in [0,T]\;\; \text{and every}\;\;
x \in \R^{d-1} \;\, \text{with some constant} \;\; K>0\;\; \text
{large enough}.
\end{align}
Assume moreover that the process $X$ can be realised as a continuous
stochastic flow $(X_s^{t,x})$ in the space variable for $s \in [0,T
\m t]$ and $(t,x) \in [0,T] \times \R^{d-1}$ and that for $z \in
\partial C$ given and fixed the following conditions are satisfied
\begin{align} \h{0pc} \label{5.12}
&\EE \Big[\;\! e^{-\Lambda_\sigma^{t,x}} G(t \p \sigma,X_\sigma^x)
\;\! \Big]  = G(t,x) + \EE \Big[ \int_0^\sigma\!\! e^{-\Lambda_s
^{t,x}} (G_t \p \LL_X G)(t \p s,X_s^x)\, ds\;\! \Big] \\ \label{5.13}
&\EE \Big[ \sup_{(t,x) \in b(z,\eps)}\, \sup_{T-t-\eps \le s \le T-t}
\, e^{-\Lambda_s^{t,x}} \big \vert \tilde H(t \p s,X_s^x) \big \vert
\;\! \Big] < \infty
\end{align}
for all stopping times $\sigma$ of $X$ with values in $[0,T \m t]$
and all $(t,x) \in b(z,\eps)$ with some $\eps
> 0$. If $z$ is probabilistically regular for $D^\circ$ then
\begin{equation} \h{6pc} \label{5.14}
\partial_t V\;\; \text{exists and is continuous at}\;\; z
\end{equation}
with $\partial_t V(z) = \partial_t G(z)$. If the hypotheses stated
above hold at every $z \in \partial C$ then $\partial_t V$ is
continuous on $[0,T]\! \times\! \R^{d-1}$.}

\vspace{12pt}

\textbf{Proof.} It will be clear from the proof below that the same
arguments are applicable in any dimension $d \ge 1$ so that for ease
of notation we will assume that $d=1$ in the sequel.

\vspace{6pt}

(I): To illustrate the arguments in a clearer manner we first
consider the special case when $\Lambda_t = 0$ for $t \ge 0$.

\vspace{6pt}

1.\ Take any $(t_n,x_n) \in C$ converging to $z \in \partial C$ as
$n \rightarrow \infty$. Passing to a subsequence of $((t_n,x_n))_{n
\ge 1}$ if needed there is no loss of generality in assuming that
\begin{equation} \h{4pc} \label{5.15}
\liminf_{n \rightarrow \infty}\, V_t(t_n,x_n) = \lim_{n \rightarrow
\infty} \frac{V(t_n \p \eps_n,x_n) \m V(t_n,x_n)}{\eps_n}
\end{equation}
for some $\eps_n \downarrow 0$ as $n \rightarrow \infty$ (we write
$V_t$ to denote $\partial V\!/\partial t$ throughout). Let $\tau_n
:= \tau_D^{t_n,x_n}$ denote the optimal stopping time for
$V(t_n,x_n)$ and set $\hat \tau_n := \tau_n \wedge (T \m t_n \m
\eps_n)$ for $n \ge 1$. Then by \eqref{5.11} and \eqref{5.12} we
find that
\begin{align} \h{0pc} \label{5.16}
&\h{2pc}V(t_n \p \eps_n,x_n) - V(t_n,x_n) \\ \notag &\h{2pc}\ge G(t_n \p \eps_n,x_n)
+ \EE \Big[ \int_0^ {\hat \tau_n}\!\! \big( G_t \p \LL_X G \p H \big)
(t_n \p \eps_n \p s,X_s^{x_n})\, ds\;\! \Big] \\ \notag &\h{2pc}\h{13pt}-
G(t_n,x_n) - \EE \Big[ \int_0^{\tau_n}\!\! \big( G_t \p \LL_X
G \p H \big)(t_n \p s,X_s^{x_n})\, ds\;\! \Big] \\ \notag &= G(t_n
\p \eps_n,x_n) - G(t_n,x_n) + \EE \Big[ \int_0^ {\hat \tau_n}\!\!
\big( \tilde H(t_n \p \eps_n \p s,X_s^{x_n}) \m \tilde H(t_n \p
s,X_s^{x_n}) \big)\, ds\;\! \Big] \\ \notag &\h{13pt}- \EE \Big[
\int_{\hat \tau_n}^{\tau_n}\!\! \tilde H(t_n \p s,X_s^{x_n})\, ds
\;\! \Big] \\[4pt] \notag &\ge G(t_n \p \eps_n,x_n) - G(t_n,x_n)
- K\, \eps_n\, \EE(\tau_n) \\[4pt] \notag &\h{13pt}- \EE \Big[\sup
_{T-t_n-\eps_n \le s \le T-t_n} \vert \tilde H(t_n \p s,X_s^{x_n})
\vert\; \eps_n\; I(T \m t_n \m \eps_n < \tau_n \le T \m t_n) \Big]
\end{align}
for $n \ge 1$. Dividing both sides by $\eps_n$ and letting $n
\rightarrow \infty$ we find from \eqref{5.15} and \eqref{5.16} that
\begin{equation} \h{7pc} \label{5.17}
\liminf_{n \rightarrow \infty}\, V_t(t_n,x_n) \ge G_t(z)
\end{equation}
where we use that $\tau_n \rightarrow 0$ almost surely as $n
\rightarrow \infty$ by probabilistic regularity of $z$ for $D^\circ$
and Corollary 6 above combined with the dominated convergence
theorem which is applicable due to \eqref{5.13} above.

\vspace{6pt}

2.\ Similarly, there is no loss of generality in assuming that
\begin{equation} \h{4pc} \label{5.18}
\limsup_{n \rightarrow \infty}\, V_t(t_n,x_n) = \lim_{n \rightarrow
\infty} \frac{V(t_n,x_n) \m V(t_n \m \eps_n,x_n)}{\eps_n}
\end{equation}
for some $\eps_n \downarrow 0$ as $n \rightarrow \infty$. By
\eqref{5.11} and \eqref{5.12} we find that
\begin{align} \h{4pc} \label{5.19}
V(t_n,x_n) - V(t_n \m \eps_n,x_n) &\le G(t_n,x_n) - G(t_n
\m \eps_n,x_n) \\ \notag &\h{13pt}\h{-11pc}+ \EE \Big[ \int_0^{\tau_n}
\big( \tilde H(t_n \p s,X_s^{x_n}) \m \tilde H(t_n \m \eps_n \p s,X_s
^{x_n}) \big)\, ds \;\! \Big] \\[4pt] \notag &\h{-11pc}\le G(t_n,x_n)
- G(t_n \m \eps_n,x_n) + K\;\! \eps_n \;\! \EE(\tau_n)
\end{align}
for $n \ge 1$. Dividing both sides by $\eps_n$ and letting $n
\rightarrow \infty$ we find from \eqref{5.18} and \eqref{5.19} that
\begin{equation} \h{7pc} \label{5.20}
\limsup_{n \rightarrow \infty}\, V_t(t_n,x_n) \le G_t(z)
\end{equation}
where we use the same arguments as following \eqref{5.17}. Combining
\eqref{5.17} and \eqref{5.20} we see that $\lim_{\,n \rightarrow
\infty} V_t(t_n,x_n) = G_t(z)$ and this completes the proof when
$\Lambda_t = 0$ for $t \ge 0$.

\vspace{6pt}

(II): Next we consider the general case when $\Lambda_t \ne 0$ for
$t \ge 0$. The proof in the general case can be carried out along
the same lines as in the special case above and we only highlight
the needed modifications throughout.

\vspace{6pt}

3.\ Taking any $(t_n,x_n) \in C$ converging to $z \in \partial C$as
$n \rightarrow \infty$ and arguing as in \eqref{5.15} above, we see
that the right-hand side of the first inequality in \eqref{5.16}
reads as follows
\begin{align} \h{1pc} \label{5.21}
&G(t_n \p \eps_n,x_n) - G(t_n,x_n) \\ \notag &\h{13pt}+ \EE \Big[
\int_0^{\hat \tau_n}\!\! \Big( e^{-\Lambda_s^{t_n+\eps_n,x_n}} \tilde
H(t_n \p \eps_n \p s,X_s^{x_n}) - e^{-\Lambda_s^{t_n,x_n}} \tilde
H(t_n \p s,X_s^{x_n}) \Big)\, ds\, \Big] \\ \notag &\h{13pt}- \EE
\Big[ \int_{\hat \tau_n}^{\tau_n}\! e^{-\Lambda_s^{t_n,x_n}} \tilde
H(t_n \p s,X_s^{x_n})\, ds\;\! \Big] \\[2pt] \notag &= G(t_n \p
\eps_n,x_n) - G(t_n,x_n) \\[2pt] \notag &\h{13pt}+ \EE \Big[ \int_0
^{\hat \tau_n}\!\! e^{-\Lambda_s^{t_n,x_n}}\! \Big( e^{\Lambda_s
^{t_n,x_n} - \Lambda_s^{t_n+\eps_n,x_n}}\! - 1 \Big)\:\! \tilde
H(t_n \p \eps_n \p s,X_s^{x_n})\, ds\, \Big] \\ \notag &\h{13pt}+
\EE \Big[ \int_0^{\hat \tau_n}\!\! \Big( e^{-\Lambda_s^{t_n,x_n}}
\Big( \tilde H(t_n \p \eps_n \p s,X_s^{x_n}) - \tilde H(t_n \p s
,X_s^{x_n}) \Big)\, ds\, \Big] \\ \notag &\h{13pt}- \EE \Big[
\int_{\hat \tau_n}^{\tau_n}\! e^{-\Lambda_s^{t_n,x_n}} \tilde
H(t_n \p s,X_s^{x_n})\, ds\;\! \Big]
\end{align}
for $n \ge 1$. The second and third expectation on the right-hand
side of \eqref{5.21} can be handled in exactly the same way as the
corresponding expectations in \eqref{5.16}, and this yields the
conclusion of \eqref{5.17}, provided that the liminf of the first
expectation on the right-hand side of \eqref{5.21} divided by
$\eps_n$ is non-negative as $n \rightarrow \infty$. To see that the
liminf is non-negative, note that \eqref{5.11} and the mean value
theorem imply that
\begin{align} \h{4pc} \label{5.22}
\frac{e^{\Lambda_s^{t_n,x_n} - \Lambda_s^{t_n+\eps_n,x_n}}-1}{\eps_n}
\ge \frac{e^{-\eps_n K s} - 1}{\eps_n} = - K s\;\! e^{-\zeta_n K s}
\end{align}
where $\zeta_n$ belongs to $(0,\eps_n)$ for $n \ge 1$. Using then
the same arguments as in \eqref{5.17} above, we see that the
inequality \eqref{5.22} yields the fact that the liminf is
non-negative so that \eqref{5.17} holds in the general case when
$\Lambda_t \ne 0$ for $t \ge 0$ as well.

\vspace{6pt}

4.\ Similarly, arguing as in \eqref{5.18} we see that the right-hand
side of the first inequality in \eqref{5.19} reads as follows
\begin{align} \h{1pc} \label{5.23}
&G(t_n,x_n) - G(t_n \m \eps_n,x_n) \\ \notag &\h{13pt}+ \EE \Big[
\int_0^{\tau_n}\!\! \Big( e^{-\Lambda_s^{t_n,x_n}} \tilde H(t_n \p
s,X_s^{x_n}) - e^{-\Lambda_s^{t_n-\eps_n,x_n}} \tilde H(t_n \m \eps_n
\p s,X_s^{x_n}) \Big)\, ds\, \Big] \\[2pt] \notag &= G(t_n,x_n) -
G(t_n \m \eps_n,x_n) \\[2pt] \notag &\h{13pt}+ \EE \Big[ \int_0
^{\tau_n}\!\! e^{-\Lambda_s^{t_n,x_n}}\! \Big( 1 - e^{\Lambda_s
^{t_n,x_n} - \Lambda_s^{t_n-\eps_n,x_n}} \Big)\:\! \tilde H(t_n
\p s,X_s^{x_n})\, ds\, \Big] \\ \notag &\h{13pt}+ \EE \Big[ \int
_0^{\tau_n}\!\! \Big( e^{-\Lambda_s^{t_n-\eps_n,x_n}} \Big( \tilde
H(t_n \p s,X_s^{x_n}) - \tilde H(t_n \m \eps_n \p s ,X_s^{x_n})
\Big)\, ds\, \Big]
\end{align}
for $n \ge 1$.

The second expectation on the right-hand side of \eqref{5.23} can be
handled in exactly the same way as the corresponding expectation in
\eqref{5.19}, and this yields the conclusion of \eqref{5.20},
provided that the limsup of the first expectation on the right-hand
side of \eqref{5.23} divided by $\eps_n$ is non-positive as $n
\rightarrow \infty$. To see that the limsup is non-positive, note
that \eqref{5.11} and the mean value theorem imply that
\begin{align} \h{4pc} \label{5.24}
\frac{1-e^{\Lambda_s^{t_n,x_n} - \Lambda_s^{t_n-\eps_n,x_n}}}{\eps_n}
\le \frac{1-e^{-\eps_n K s}}{\eps_n} = K s\;\! e^{-\zeta_n K s}
\end{align}
where $\zeta_n$ belongs to $(0,\eps_n)$ for $n \ge 1$. Using then
the same arguments as in \eqref{5.20} above, we see that the
inequality \eqref{5.24} yields the fact that the limsup is
non-positive so that \eqref{5.20} holds in the general case when
$\Lambda_t \ne 0$ for $t \ge 0$ as well. Combining the conclusions
of \eqref{5.17} and \eqref{5.20} we see that $\lim_{\,n \rightarrow
\infty} V_t(t_n,x_n) = G_t(z)$ and this completes the proof. \hfill
$\square$

\vspace{16pt}

\textbf{Remark 16.} Note that the comments on the sufficient
conditions from Theorem 8 made in Remark 9 above extend to the
corresponding sufficient conditions in Theorem 15 and we will omit
further details in this direction. Note also that the proof of
\eqref{5.20} above could also be accomplished by means of the mean
value theorem (as in the proof of Theorems 8 and 10) without appeal
to the identity \eqref{5.12}.

\vspace{12pt}

3.\ The result and proof of Theorem 13 and Theorem 15 extend to the
case when the gain function $G$ in the optimal stopping problem
\eqref{2.1}/\eqref{2.2} is \emph{not smooth} away from the optimal
stopping boundary $\partial C$. Instead of formulating a general
theorem of this kind, which would be overly technical and rather
difficult to read, we will illustrate key arguments of such
extensions through an important example that was already considered
in Example 12 above for the space derivative. The method of proof to
be presented below is different from the method of proof applied in
Example 12 above.

\vspace{12pt}

\textbf{Example 17 (Continuity of the time derivative in the
American put)}. Consider the optimal stopping problem \eqref{4.28}
above where $X$ solves \eqref{4.29}. Horizon in the optimal stopping
problem \eqref{4.28} is finite so that the setting belongs to
Theorem 15 above. Since the gain function $G(x) = (K \m x)^+$ for
$x>0$ is not differentiable at $K$ we see that the condition
\eqref{5.10} fails and hence we cannot conclude that
\begin{equation} \h{7pc} \label{5.25}
V_t\;\; \text{is continuous on}\;\; \partial C
\end{equation}
using Theorem 15 (we write $V_t$ to denote $\partial V\!/\partial t$
throughout). We will now show however that the method of proof of
Theorem 15 extends to cover the case of the non-differentiable gain
function $G(x) = (K \m x)^+$ for $x>0$. This will also serve as an
illustration of how similar other cases of non-smooth gain functions
$G$ in the optimal stopping problem \eqref{2.1}/\eqref{2.2} can be
handled. The derivation of \eqref{5.25} will be divided in three
steps as follows.

\vspace{6pt}

1.\ We first recall the facts about the optimal stopping problem
\eqref{4.28} stated in the first step of the proof of \eqref{4.30}
above. In particular, taking any sequence $(t_n,x_n) \in C$
converging to $z=(t,b(t)) \in \partial C$ we know that
$\tau_D^{t_n,x_n} \rightarrow 0$ almost surely as $n \rightarrow
\infty$. Moreover, applying the It\^o-Tanaka formula, we find using
\eqref{4.29} that
\begin{align} \h{0pc} \label{5.26}
e^{-rt} (K \m X_t)^+ &= (K \m x)^+ - \int_0^t r e^{-rs}\:\! K\:\!
I(X_s\! <\! K)\, ds - \int_0^t e^{-rs}\:\! \sigma\:\! X_s\;\! I(
X_s\! <\! K)\, dB_s \\ \notag &\h{13pt}+ \int_0^t \tfrac{1}{2}
e^{-rs}\, d \ell_s^K(X)
\end{align}
for $t \in [0,T]$ where $\ell^K(X)$ is the local time process of $X$
defined by
\begin{equation} \h{3pc} \label{5.27}
\ell_t^K(X) = \lim_{\eps \downarrow 0}\, \frac{1}{2 \eps} \int_0^t
I(K \m \eps < X_s < K \p \eps)\: d \langle X,X \rangle_s
\end{equation}
where the convergence takes place in probability and the quadratic
variation process $\langle X,X \rangle$ of $X$ is given by $\langle
X,X \rangle_t = \int_0^t \sigma^2 X_s^2\, ds$ for $t \in [0,T]$. It
is easily verified that the third term on the right-hand side in
\eqref{5.26} defines a continuous martingale for $t \in [0,T]$.
Hence by the optional sampling theorem we find that the Bolza
formulated optimal stopping problem \eqref{4.28} can be Lagrange
reformulated (see \cite[p.\ 141]{PS} for the terminology) as follows
\begin{align} \h{2pc} \label{5.28}
\tilde V(t,x) &:= V(t,x) - (K \m x)^+ \\ \notag &= \sup_{0 \le \tau
\le T \m t} \EE \Big[ \int_0^\tau \tfrac{1}{2} e^{-rs}\, d\ell_s^K
(X^x) - \int_0^ \tau r e^{-rs}\:\! K\:\! I(X_s^x\! <\! K) \, ds
\;\! \Big]
\end{align}
for $t \in [0,T]$ and $x>0$. Thus the optimal stopping problems
\eqref{4.28} and \eqref{5.28} are equivalent and a stopping time is
optimal in \eqref{4.28} if and only if it is optimal in
\eqref{5.28}.

\vspace{6pt}

2.\ We next connect to the first part of the proof of Theorem 15.
Passing to a subsequence of $((t_n,x_n))_{n \ge 1}$ if needed there
is no loss of generality in assuming that
\begin{equation} \h{4pc} \label{5.29}
\liminf_{n \rightarrow \infty}\, V_t(t_n,x_n) = \lim_{n \rightarrow
\infty} \frac{V(t_n \p \eps_n,x_n) \m V(t_n,x_n)}{\eps_n}
\end{equation}
for some $\eps_n \downarrow 0$ as $n \rightarrow \infty$. Let
$\tau_n := \tau_D^{t_n,x_n}$ be the optimal stopping time for
$V(t_n,x_n)$ and thus $\tilde V(t_n,x_n)$ as well. Set $\hat \tau_n
:= \tau_n \wedge (T \m t_n \m \eps_n)$ for $n \ge 1$. We then have
\begin{align} \h{0pc} \label{5.30}
V(t_n \p \eps_n,x_n) - V(t_n,x_n) &= \tilde V(t_n \p \eps_n,x_n) -
\tilde V(t_n,x_n) \\[2pt] \notag &\ge \EE \Big[ \int_0^{\hat \tau_n}
\tfrac{1}{2} e^{-rs}\, d\ell_s^K(X^{x_n}) - \int_0^{\hat \tau_n} r e^
{-rs}\:\! K\:\! I(X_s^{x_n}\! <\! K)\, ds\;\! \Big] \\ \notag &\h{13pt}
- \EE \Big[ \int_0^{\tau_n} \tfrac{1}{2} e^{-rs}\, d\ell_s^K(X^{x_n})
- \int_0^{\tau_n} r e^{-rs}\:\! K\:\! I(X_s^{x_n}\! <\! K)\, ds\;\!
\Big] \\ \notag &\ge -\EE \Big[ \int_{\hat \tau_n}^{\tau_n} \tfrac{1}
{2} e^{-rs}\, d \ell_s^K(X^{x_n})\; I(T \m t_n \m \eps_n < \tau_n \le
T \m t_n)\, \Big] \\[3pt] \notag &\ge -\frac{1}{2}\, e^{-r(T-t_n-
\eps_n)}\, \EE \big[ \ell_{T-t_n}^K(X^{x_n}) - \ell_{T-t_n-\eps_n}
^K(X^{x_n}) \big]
\end{align}
for all $n \ge 1$. By \eqref{5.27} and Fatou's lemma we find that
\begin{align} \h{2pc} \label{5.31}
&\EE \big[ \ell_{T-t_n}^K(X^{x_n}) - \ell_{T-t_n-\eps_n}^K(X
^{x_n}) \big] \\ \notag &\h{3pc}= \EE \Big[ \lim_{\eps \downarrow 0}\,
\frac{1}{2 \eps} \int_{T-t_n-\eps_n}^{T-t_n} I(K \m \eps < X_s^{x_n}
< K \p \eps)\: \sigma^2 (X_s^{x_n})^2\, ds\;\! \Big] \\ \notag
&\h{3pc} \le \sigma^2 x_n^2\, \liminf_{\eps \downarrow 0} \int_{T-t
_n-\eps_n}^{T-t_n} \frac{1}{2 \eps}\, \EE \big[\, I \big( \tfrac{K
-\eps}{x_n} < X_s^1 < \tfrac{K + \eps}{x_n} \big)\:\! (X_s^1)^2\;\!
\big]\, ds \\ \notag &\h{3pc}= \sigma^2 x_n^2\, \liminf_{\eps
\downarrow 0} \int_{T-t_n-\eps_n}^{T-t_n} \Big(\;\! \frac{1}{2 \eps}
\int_{\frac{K - \eps}{x_n}}^{\frac{K + \eps}{x_n}} x^2 f_{X_s^1}(x)\,
dx \Big)\, ds \\ \notag &\h{3pc}= \sigma^2 K^2\! \int_{T-t_n-\eps_n}
^{T-t_n}\! f_{X_s^1}(\tfrac{K}{x_n})\, ds
\end{align}
for all $n \ge 1$ where $f_{X_s^1}$ denotes the density function of
$X_s^1$ for $s>0$ and in the last equality we use the dominated
convergence theorem. Using the scaling property $B_s \sim
\sqrt{s}\;\! B_1$ it is easily verified that $f_{X_s^1}$ is given by
\begin{equation} \h{5pc} \label{5.32}
f_{X_s^1}(x) = \frac{1}{\sigma x \sqrt{s}}\: \varphi \bigg( \frac{
\log(x) \m (r \m \sigma^2/2)s}{\sigma \sqrt{s}} \bigg)
\end{equation}
for $x>0$ and $s>0$ where $\varphi$ denotes the standard normal
density function given by $\varphi(x) = (1/\sqrt{2 \pi})\,
e^{-x^2/2}$ for $x \in \R$. Inserting \eqref{5.32} into \eqref{5.31}
we find that
\begin{equation} \h{5pc} \label{5.33}
\EE \big[ \ell_{T-t_n}^K(X^{x_n}) - \ell_{T-t_n-\eps_n}^K(X^{x_n})
\big] \le c\, \eps_n
\end{equation}
for all $n \ge n_0$ with some $n_0 \ge 1$ large enough, where the
constant $c = c(T \m t)$ is given by
\begin{equation} \h{5pc} \label{5.34}
c = \sigma K^2\, \sup\: \frac{e^{-y}}{\sqrt{s}}\, \varphi \bigg(
\frac{ y \m (r \m \sigma^2/2)s}{\sigma \sqrt{s}} \bigg)
\end{equation}
with the supremum being taken over all $s \in [(T-t)/2,2(T-t)]$ and
$y \in \R$ (upon substituting $y = \log(x)$ in \eqref{5.32} above).
Making use of \eqref{5.33} in \eqref{5.31} we obtain
\begin{align} \h{5pc} \label{5.35}
V(t_n \p \eps_n,x_n) - V(t_n,x_n) \ge -c\, \eps_n
\end{align}
for all $n \ge n_0$. Note that we can formally replace $x_n$ in
\eqref{5.35} by $x$ because the constant $c$ depends only on $T \m t
> 0$ and the resulting inequality holds uniformly over all $x > 0$.

Having \eqref{5.35} we modify the optimal stopping time $\tau_n$ by
setting $\tau_n^\delta := \tau_n \wedge \delta$ where $\delta > 0$
is any (small) number such that $t_n \p \eps_n \p \delta \le T$ for
all $n \ge n_1$ where $n_1 \ge n_0$ is sufficiently large. (Note
that this is possible since $t<T$ with $t_n \rightarrow t$ and
$\eps_n \downarrow 0$ as $n \rightarrow \infty$.) Since $(t,x)
\mapsto e^{-rt} V(t,x)$ is superharmonic on $[0,T]\! \times\!
(0,\infty)$ and harmonic on $C$, we find that
\begin{align} \h{1.5pc} \label{5.36}
&V(t_n \p \eps_n,x_n) - V(t_n,x_n) \ge \EE \Big[ e^{-r \tau_n^\delta}
\Big( V(t_n \p \eps_n \p \tau_n^\delta,X_{\tau_n^\delta}^{x_n}) -
V(t_n \p \tau_n^\delta,X_{\tau_n^\delta}^{x_n}) \Big) \Big] \\[1pt]
\notag &= \EE \Big[ e^{-r \tau_n} \Big( V(t_n \p \eps_n \p \tau_n
,X_{\tau_n}^{x_n}) - (K \m X_{\tau_n}^{x_n})^+) I( \tau_n \le
\delta) \Big] \\[1pt] \notag &\h{13pt}+ \EE \Big[ e^{-r \delta}
\Big( V(t_n \p \eps_n \p \delta ,X_{\delta}^{x_n}) - V(t_n \p
\delta,X_{\delta}^{x_n}) \Big) I( \tau_n > \delta) \Big] \\[1pt]
\notag &\ge -c\, \eps_n\:\! \PP( \tau_n > \delta)
\end{align}
for all $n \ge n_1$ where in the final inequality we use
\eqref{5.35} applied to $(t_n \p \delta,x)$ in place of $(t_n,x_n)$
for $n \ge 1$ and holding uniformly over all $x>0$. Dividing both
sides in \eqref{5.36} by $\eps_n$ we find from \eqref{5.29} that
\begin{equation} \h{5pc} \label{5.37}
\liminf_{n \rightarrow \infty}\, V_t(t_n,x_n) \ge 0 = G_t(z)
\end{equation}
where we use that $\tau_n \rightarrow 0$ almost surely so that $\PP(
\tau_n > \delta) \rightarrow 0$ as $n \rightarrow \infty$.

\vspace{6pt}

3.\ We finally connect to the second part of the proof of Theorem
15. Similarly, there is no loss of generality in assuming that
\begin{equation} \h{3pc} \label{5.38}
\limsup_{n \rightarrow \infty}\, V_t(t_n,x_n) = \lim_{n \rightarrow
\infty} \frac{V(t_n,x_n) \m V(t_n \m \eps_n,x_n)}{\eps_n}
\end{equation}
for some $\eps_n \downarrow 0$ as $n \rightarrow \infty$. We then
have
\begin{align} \h{0pc} \label{5.39}
V(t_n,x_n) - V(t_n \m \eps_n,x_n) &= \tilde V(t_n,x_n) - \tilde
V(t_n \m \eps_n,x_n) \\[2pt] \notag &\le \EE \Big[ \int_0^{\tau_n}
\tfrac{1}{2} e^{-rs}\, d\ell_s^K(X^{x_n}) - \int_0^{\tau_n} r e
^{-rs}\:\! K\:\! I(X_s^{x_n}\! <\! K)\, ds\;\!
\Big] \\ \notag &\h{13pt}- \EE \Big[ \int_0^{\tau_n} \tfrac{1}{2}
e^{-rs}\, d\ell_s^K(X^{x_n}) - \int_0^{\tau_n} r e^{-rs}\:\! K\:\!
I(X_s^{x_n}\! <\! K)\, ds\;\! \Big] \\ \notag &= 0
\end{align}
for $n \ge 1$. Note that this inequality also follows from
\eqref{4.28} from where we see directly that $t \mapsto V(t,x)$ is
decreasing on $[0,T]$ for $x>0$. Dividing both sides in \eqref{5.39}
by $\eps_n$ we find from \eqref{5.38} and \eqref{5.39} that
\begin{equation} \h{6pc} \label{5.40}
\limsup_{n \rightarrow \infty}\, V_t(t_n,x_n) \le 0 = G_t(z)\, .
\end{equation}
Combining \eqref{5.37} and \eqref{5.40} we see that $\lim_{\,n
\rightarrow \infty} V_t(t_n,x_n) = 0 = G_t(z)$ so that \eqref{5.25}
holds as claimed and the proof is complete. \hfill $\square$

\vspace{12pt}

\textbf{Remark 18.} Note that the method of proof presented in
Example 17 first derives Lipschitz continuity of $t \mapsto V(t,x)$
uniformly over all $x$ and then `lifts' this continuity to $C^1$
regularity of $t \mapsto V(t,x)$ at $z \in \partial C$ using the
superharmonic property of $(t,x) \mapsto e^{-rt} V(t,x)$ on $[0,T]\!
\times\! (0,\infty)$. To our knowledge this `lifting' method is
applied in Example 17 for the first time in the literature. In
addition to yielding the first known probabilistic proof of
\eqref{5.25} in the American put problem, it is also clear from the
arguments used in Example 17 that the `lifting' method is applicable
to a large class of diffusion/Markov processes in optimal stopping
and free boundary problems with non-smooth gain functions.

\vspace{16pt}

\textbf{Acknowledgements.} T.\ De Angelis gratefully acknowledges
partial support by EPSRC Grant EP/R021201/1 while working on the
paper.

\vspace{6pt}


\begin{center}

\end{center}


\par \leftskip=24pt

\vspace{20pt}

\nt Tiziano De Angelis \\
School of Mathematics \\
University of Leeds \\
Leeds LS2 9JT \\
United Kingdom \\
\texttt{t.deangelis@leeds.ac.uk}

\leftskip=25pc \vspace{-87pt}

\nt Goran Peskir \\
Department of Mathematics \\
The University of Manchester \\
Oxford Road \\
Manchester M13 9PL \\
United Kingdom \\
\texttt{goran@maths.man.ac.uk}

\par

\end{document}